\documentclass{siamart220329}

\usepackage{graphicx}
\usepackage{xcolor}
\usepackage{textcomp}
\usepackage{booktabs}
\usepackage{multirow}
\usepackage{mathrsfs}
\usepackage{algorithm}
\usepackage{algorithmic}   
\usepackage{listings}
\usepackage{subfig}       
\usepackage{float}           
\usepackage{placeins}
\usepackage{adjustbox}
\usepackage{siunitx}
\usepackage{pdfpages}
\usepackage{pdflscape}
\usepackage{tabularx}
\usepackage{amssymb}

\title{Efficient Multi-Precision Computation of Bessel Functions for Real Orders and Complex Arguments with a Fortran Implementation---Part~III: Regular Bessel Functions of the First and Second Kinds \texorpdfstring{\(J_{\nu}(z)\) and \(Y_{\nu}(z)\)}{Jν(z) and Yν(z)}}

\author{
  Mofreh R.~Zaghloul%
  \thanks{Department of Physics, College of Sciences,
    United Arab Emirates University, Al~Ain 15551, Abu Dhabi, UAE
    (\email{m.zaghloul@uaeu.ac.ae}).}
  \thanks{Department of Mathematics,
    Massachusetts Institute of Technology, Cambridge, MA, USA.}
  \and
  Steven G.~Johnson%
  \thanks{Department of Mathematics,
    Massachusetts Institute of Technology,
    182 Memorial Dr., Cambridge, MA 02139, USA
    (\email{stevenj@math.mit.edu}).}
}

\headers{Efficient Computation of Bessel Functions}{M. R.~Zaghloul and S. G.~Johnson}

\ifpdf
\hypersetup{
  pdftitle={Efficient and Accurate Computation of Bessel Functions},
  pdfauthor={Mofreh R.~Zaghloul and Steven G.~Johnson}
}
\fi

\newcommand{\Jnu}{J_{\nu}}
\newcommand{\Ynu}{Y_{\nu}}

\begin{document}

\maketitle


\begin{abstract}
This paper is the third and final part in a series devoted to the development of numerically stable and efficient algorithms, together with multi-precision
\texttt{Fortran} implementations, for computing Bessel functions with real orders and complex arguments. Parts~I and~II presented stable and efficient algorithms for the modified Bessel functions $I_{\nu}(z)$ and $K_{\nu}(z)$; here we treat the regular Bessel functions $\Jnu(z)$ and $\Ynu(z)$.

The proposed algorithms support complex arguments for both positive and \textit{negative} real orders and are implemented in native double and quadruple precision. Quadruple (quad) precision substantially increases dynamic range and accuracy (by roughly an order of magnitude in reliably computable \(|\nu|\) and \(|z|\)), thereby extending applicability to problems requiring 20--30 significant digits.

Comprehensive accuracy and performance comparisons are carried out against both the widely used AMOS package (Algorithm~644), which is restricted to double-precision arithmetic, and the more recent Algorithm~912, which supports both double- and quadruple-precision arithmetic as well as complex orders. In double precision, the present implementation consistently outperforms Algorithm~644, achieving execution times of approximately $35\!-\!67\%$ for $\Jnu(z)$ and $44\!-\!72\%$ for $\Ynu(z)$, while also producing reliable results in parameter regimes where Algorithm~644 fails (e.g., negative orders).

In comparison with Algorithm~912, the present algorithm achieves comparable accuracy in double-precision computations and significantly higher accuracy in quadruple-precision calculations. At the same time, it requires only a small fraction of the computational cost, typically ranging from a few thousandths to a few hundredths of the time taken by Algorithm~912, depending on the precision and the parameter regime. Furthermore, unlike Algorithm~912, whose applicability is restricted to a limited region of the \((\Re \nu, z)\)-plane, the present algorithm remains stable and accurate over the full tested domain for which reliable reference values are available.

The implementation and associated test data will be deposited in an open-access repository to support reproducible research.
\end{abstract}



\section{Introduction}\label{sec:intro}

The general Bessel differential equation,
\begin{equation}
z^2 \frac{d^2 y}{dz^2} + z \frac{dy}{dz} + (z^2 - \nu^2)y = 0,
\label{eq:Bessel_diff_eq}
\end{equation}
admits two linearly independent solutions known as the \emph{Bessel functions} of the first and second kinds, \(\Jnu(z)\) and \(\Ynu(z)\) \cite{Abramowitz1964,Oldham2009,Olver2010,Temme1996,Watson1922}.  For complex arguments \(z = x + \mathrm{i}y\), these functions extend naturally to the Hankel functions \cite{NIST2025,Oldham2009,Olver2010}
\begin{equation}
H_{\nu}^{(1)}(z) = \Jnu(z) + \mathrm{i}\,\Ynu(z), \qquad
H_{\nu}^{(2)}(z) = \Jnu(z) - \mathrm{i}\,\Ynu(z),
\label{eq:Hankel_def}
\end{equation}
which represent outgoing and incoming cylindrical waves. 
\noindent
Together, \(\Jnu\), \(\Ynu\), \(H_{\nu}^{(1)}\), and \(H_{\nu}^{(2)}\) form a fundamental basis for describing oscillatory and radiating wave 
phenomena in applied mathematics, physics, and engineering \cite{Zhang1996,NIST2025,Thompson1986}.

\medskip
\noindent
The classical literature provides extensive numerical methods for Bessel functions, particularly for real arguments and moderate parameter ranges 
\cite{Abramowitz1964,Gautschi1964,Amos1973,Amos1974,Amos1974b,Amos1977,Temme1973,Boisvert1992,Segura2011,Cody1993}. However, comparatively few studies address the accurate and robust computation 
of \(\Jnu(z)\) and \(\Ynu(z)\) for real~\(\nu\) and genuinely complex~\(z\), especially in parameter regimes requiring high accuracy; the most widely used of these is the complex-argument package of Amos \cite{Amos1986}. In addition to Algorithm~644, more recent work has addressed the computation of Bessel functions with enhanced precision and generality. In particular, Algorithm~912 \cite{Kodama2011} provides routines for Bessel functions of complex order and complex argument, with support for both double- and quadruple-precision arithmetic. However, its practical applicability is limited in parameter range, and its computational cost can be relatively high, which restricts its use in large-scale or high-performance applications.
 
\medskip
\noindent
Complex-plane computation presents distinctive challenges stemming from oscillatory behavior, exponential growth and decay in different quadrants, and 
the branch structure introduced by analytic continuation \cite{Thompson1986,GilSeguraTemme2007}.  

\medskip
\noindent
As noted in the first two parts of this series \cite{Zaghloul_Johnson_I,Zaghloul_Johnson_II}, Amos’s Algorithm~644 \cite{Amos1986}, which underlies many standard numerical libraries, is limited to double-precision arithmetic. In addition, it does not support negative real orders and frequently triggers premature overflow or underflow in the complex plane. Its loss of accuracy for large \(|\nu|\) or \(|z|\), together with the absence of quadruple-precision support, significantly restricts its applicability in many of the
modern high-accuracy simulations in electromagnetics, acoustics, photonics, quantum scattering, and plasma modeling, where high accuracy over wide parameter
ranges is essential. Such applications increasingly require 20--30 accurate digits and robust treatment of large, oscillatory parameter regimes. 

\medskip
\noindent
Quadruple-precision support in a compiled-language implementation provides a practical compromise between the limited range, limited accuracy of double precision and the
high computational cost of arbitrary-precision environments (e.g. \texttt{Maple} \cite{MAPLE2015} or \texttt{mpmath}) \cite{mpmath}. To the best of the authors' knowledge, no publicly available compiled-language implementation currently provides native quadruple-precision evaluation of \(\Jnu(z)\) and \(\Ynu(z)\) for complex arguments and real orders, including negative orders, over comparably broad tested domains.

\medskip
\noindent
To address these limitations, we introduce new algorithms for these functions designed to deliver high accuracy, numerical stability, and computational efficiency over the tested double- and quadruple-precision domains.

\medskip
\noindent
This paper constitutes the third and final part of a unified series devoted to the accurate and efficient computation of Bessel and modified Bessel functions of real positive and negative orders with complex arguments, using native multiprecision arithmetic. Part~I \cite{Zaghloul_Johnson_I} developed a robust algorithm for the modified 
Bessel function \(I_{\nu}(z)\), and Part~II \cite{Zaghloul_Johnson_II} extended the methodology to \(K_{\nu}(z)\). The present work completes this framework by treating the regular Bessel functions of the first and second kinds, \(\Jnu(z)\) and \(\Ynu(z)\).

\medskip
\noindent
These functions are related to the modified Bessel functions through analytic continuation \cite{Temme1996,Olver2010}:
\begin{equation}
\Jnu(z) = e^{+\mathrm{i}\pi\nu/2}\,I_{\nu}(-\mathrm{i}z), \qquad
\Ynu(z) = \frac{\cos(\pi\nu)\,\Jnu(z) - J_{-\nu}(z)}{\sin(\pi\nu)},
\qquad \nu \notin \mathbb{Z},
\label{eq:JY_relations_intro}
\end{equation}

\noindent
which provide an analytic bridge between the modified and regular Bessel functions.  

\medskip
\noindent
Direct numerical use of these relations is nontrivial: exponential phase factors can amplify round-off error, and the difference formulas for \(\Ynu(z)\) 
become ill-conditioned near integer orders \cite{Temme1973,GilSeguraTemme2007}.  
To overcome these difficulties, we develop hybrid methods that combine stabilized analytic continuation with region-specific numerical representations.

\medskip
\noindent
The approach integrates Frobenius-type power series for small arguments, uniform asymptotic (Debye-type) expansions for large orders or arguments 
\cite{FrenzenWong1984,Temme1973}, and forward/backward recurrence relations for intermediate regimes \cite{GilSeguraTemme2011}.  
The algorithms automatically select the most stable representation based on \(|z|\) and \(|\nu|\), with special care to maintain continuous phase tracking 
across quadrants \cite{Thompson1986}. The resulting framework yields accurate, stable evaluation of \(\Jnu(z)\) and \(\Ynu(z)\) over the tested complex-domain grids.

\medskip
\noindent
As in the preceding parts of this series, the algorithms are implemented in modern, modular \texttt{Fortran} and support both double and quadruple precision. Validation against high-precision \texttt{Maple} and \texttt{mpmath} data \cite{MAPLE2015,mpmath} confirms stability and accuracy across the entire tested domain. In double precision the new algorithms exhibit improved numerical robustness and (typically) reduced wall-clock time compared with Algorithm~644 on our test platforms; in quadruple precision, the present methods extend the reliable computational domain substantially (relative to double precision), as Algorithm~644 does not support quadruple precision.

\section{Algorithms}\label{sec:algorithm}
This section describes the numerical schemes employed for the evaluation of the Bessel
functions of real order and complex argument.  The algorithms combine power–series,
asymptotic expansions, and analytic continuation in a manner that ensures numerical
stability throughout the complex plane.  Each function is treated in its own subsection.

\subsection{Bessel Functions of the First Kind, $J_\nu(z)$}\label{sec:J_nu}
The solution of~\eqref{eq:Bessel_diff_eq} is analytic in \(z\in\mathbb{C}\) except for a branch point at \(z=0\) when \(\nu\notin\mathbb{Z}\).
On the principal branch, \(\Jnu(z)\) is defined via the principal value of \((z/2)^{\nu}\) and is analytic in the complex \(z\)-plane cut along the negative real axis, \((-\infty,0]\)~\cite{NIST2025}.
For integer orders, \(\Jnu(z)\) is entire in \(z\); for noninteger orders, the branch cut ensures the single-valuedness of \((z/2)^{\nu}\).
Moreover, for \(\nu\notin\mathbb{Z}\) the pair \(\{\Jnu,J_{-\nu}\}\) forms two linearly independent solutions.

\medskip
\noindent
Applying the Frobenius method at \(z=0\) yields two power-series solutions with leading exponents \(\pm\nu\), corresponding to \(\Jnu(z)\) and \(J_{-\nu}(z)\).
As in the case of \(I_{\nu}(z)\)~\cite{Zaghloul_Johnson_I}, the resulting expansion provides an exact analytic representation that converges for all finite \(z\) and is especially efficient for small arguments~\cite{Abramowitz1964,Olver2010,Luke1969}.
For either choice of sign, the series reads
\begin{equation}
\label{eq:J_series_pm}
J_{\pm\nu}(z)
= \left(\frac{z}{2}\right)^{\pm\nu}
  \sum_{k=0}^{\infty}
  \frac{(-1)^{k}}{k!\,\Gamma(k\!\pm\!\nu+1)}
  \left(\frac{z}{2}\right)^{2k}.
\end{equation}

\noindent
This series is analogous to the small-\(z\) expansion used for evaluating \(I_{\nu}(z)\) in~\cite{Zaghloul_Johnson_I}, differing only by the alternating factor \((-1)^{k}\).
It can therefore be accumulated through a similar stable recurrence:
\begin{equation}
\label{eq:J_series_recurrence}
T_{0}=1,\qquad
T_{k+1}=-\,\frac{(z/2)^{2}}{(k+1)(k\!\pm\!\nu+1)}\,T_{k},
\qquad k = 0, 1, 2, \ldots
\end{equation}
and
\begin{equation}
\label{eq:J_series_sum}
J_{\pm\nu}(z)
= \frac{\bigl(z/2\bigr)^{\pm\nu}}{\Gamma(1\!\pm\!\nu)}
  \sum_{k=0}^{\infty} T_{k}.
\end{equation}

\noindent
Consistent with Part~I~\cite{Zaghloul_Johnson_I}, the convergence of this series is absolute for all finite $z$; however, the \emph{rate} of convergence depends on both $|z|$ and $|\nu|$. Owing to the structural similarity between the power-series representations of $\Jnu(z)$ and $I_{\nu}(z)$, which differ essentially only by an alternating sign, the same underlying routine, \texttt{bessel\_series\_core}, is employed for both functions, with a simple integer flag selecting the appropriate form. For $|z|\ll1$, only a few terms are required to achieve high accuracy, whereas for larger $|z|$ the alternating terms decrease more slowly, and numerical cancellation may arise. Consequently, and analogous to the treatment of $I_{\nu}(z)$, this expansion is used only in the small-argument region, typically defined as $|z|^{2}\lesssim 324 + 8|\nu|$, where it provides an optimal balance of accuracy and efficiency. Together with the analytic continuation described below, the same formulation applies to both positive and negative orders over the tested complex-domain grids.

\medskip
\noindent
For the remaining region of the \(\nu\)–\(|z|\) domain, we employ the numerically stable algorithm and routine for \(I_{\nu}(z)\) developed in~\cite{Zaghloul_Johnson_I}.
For all real orders \(\nu\), positive and negative, and complex arguments \(z\), \(\Jnu(z)\) can be obtained from \(I_{\nu}(z)\) through the analytic-continuation identity
\begin{equation}\label{eq:JtoI}
  J_\nu(z)=
  \begin{cases}
    \exp\!\big(+\mathrm{i}\tfrac{\pi\nu}{2}\big)\, I_\nu(-\mathrm{i} z), & \operatorname{Im}(z)\ge 0,\\[3pt]
    \exp\!\big(-\mathrm{i}\tfrac{\pi\nu}{2}\big)\, I_\nu(+\mathrm{i} z), & \operatorname{Im}(z)< 0,
  \end{cases}
\end{equation}
which is valid for all real \(\nu\) under the principal branch of \(\log z\).
The phase factor $\exp(\pm i\pi\nu/2)$ naturally accounts for both positive and negative orders, requiring no special handling in the implementation. A pseudo-algorithm for the method is provided in Algorithm~\ref{alg:Jnu-series-mapping}.

\begin{algorithm}[H]
\caption{Evaluation of $J_{\nu}(z)$ for real $\nu$ and complex $z$
(series for small $|z|$; mapping via $I_{\nu}$ otherwise, using
\texttt{i\_nu\_of\_z} from Part~I \cite{Zaghloul_Johnson_I})}
\label{alg:Jnu-series-mapping}
\begin{algorithmic}[1]
  \STATE \textbf{Input:} real order $\nu$, complex argument $z$.
  \STATE \textbf{Output:} $J_{\nu}(z)$ in \texttt{J}, status \texttt{ierr}.
  \STATE Set $x\gets\Re(z)$, $y\gets\Im(z)$, $r^2\gets x^2+y^2$, $\nu_a\gets|\nu|$, \texttt{ierr}$\gets0$.
  \STATE Set $n\gets\operatorname{NINT}(\nu)$, $f\gets\nu-n$ \quad (so $f\in[-\tfrac12,\tfrac12]$).
  \STATE Set $\texttt{small\_z\_border}\gets 324+8\nu_a$.
  \STATE Set $(z_I,y_{\text{sign}})\gets(-\mathrm{i}z,+1)$ if $y\ge0$, otherwise set $(z_I,y_{\text{sign}})\gets(+\mathrm{i}z,-1)$.

  \IF{$(\nu\ge0)$ and $(r^2\le\texttt{small\_z\_border})$}
    \STATE Set $\texttt{phase}\gets1$, $z_{\text{series}}\gets z$.
    \STATE Call \texttt{bessel\_series\_core}$(\nu,z_{\text{series}},-1)\to(\texttt{J},\texttt{ierr})$.
    \STATE Set $J_{\nu}(z)\gets\texttt{phase}\cdot\texttt{J}$ and return.
  \ELSE
    \STATE Compute $m_4\gets\operatorname{modulo}(n,4)$ and set
    \[
    \texttt{base}\gets
    \begin{cases}
      1, & m_4=0,\\
      \mathrm{i}, & m_4=1,\\
      -1, & m_4=2,\\
      -\mathrm{i}, & m_4=3.
    \end{cases}
    \]
    \STATE If $y_{\text{sign}}<0$, set $\texttt{base}\gets\overline{\texttt{base}}$.
    \STATE Set $\phi\gets y_{\text{sign}}\,f\,\pi/2$.
    \STATE Set $\texttt{frac\_phase}\gets\cos\phi+\mathrm{i}\sin\phi$.
    \STATE Set $\texttt{csgn}\gets\texttt{base}\cdot\texttt{frac\_phase}$.

    \IF{$\Im(z_I)<0$}
      \STATE Set $\widetilde z_I\gets\overline{z_I}$.
      \STATE Compute $\widetilde I\gets\texttt{i\_nu\_of\_z}(\nu,\widetilde z_I;\texttt{ierr})$.
      \STATE Set $I\gets\overline{\widetilde I}$.
    \ELSE
      \STATE Compute $I\gets\texttt{i\_nu\_of\_z}(\nu,z_I;\texttt{ierr})$.
    \ENDIF

    \STATE If \texttt{ierr}$\ne0$, return.
    \STATE Set $J_{\nu}(z)\gets\texttt{csgn}\cdot I$ and return.
  \ENDIF
\end{algorithmic}
\end{algorithm}

\noindent
The split $\nu=n+f$ ensures the continuous rotation is small ($|f|\le \tfrac{1}{2}$) and the discontinuous part is exact from the $n\bmod 4$ table.

\subsection{Bessel Function of the Second Kind, $Y_\nu(z)$}\label{sec:Y_nu}

The Bessel function of the second kind, $Y_\nu(z)$, is evaluated for real orders $\nu\in\mathbb{R}$ and complex arguments $z\in\mathbb{C}$ using a hybrid strategy that combines analytic reflection identities with a branch-corrected mapping through the modified Bessel functions. This approach provides numerical stability across the entire complex plane and avoids the singularities associated with direct formulas containing $\sin(\pi\nu)$ in the denominator.

\medskip
\noindent
The function $\Ynu(z)$ can be expressed in terms of the Bessel functions of the first kind, $\Jnu(z)$ and $J_{-\nu}(z)$, through the reflection relation (second part of~\eqref{eq:JY_relations_intro}), that is,

\begin{equation}\label{eq:Y_reflection}
\Ynu(z)
= \frac{\cos(\pi\nu)\,\Jnu(z) - J_{-\nu}(z)}{\sin(\pi\nu)},
\qquad \nu \notin \mathbb{Z},
\end{equation}
which is numerically well conditioned except near integer orders, where \(\sin(\pi\nu)\to 0\).
In those regions, we employ an analytic continuation through the modified Bessel functions using the relation
\begin{equation}\label{eq:Y_K_mapping}
\Ynu(z)
= -{\mathrm{i}}\left[H^{(1)}_{\nu}(z) - \Jnu(z)\right].
\end{equation}
\noindent
where, as introduced in the Introduction, \(H^{(1)}_{\nu}(z)\) is the Hankel function of the first kind. The Hankel functions of the first and second kinds are related to the modified Bessel functions of the second kind by
\begin{equation}
H_{\nu}^{(1)}(z) = \frac{2}{\pi\,\mathrm{i}}\,
  e^{-\mathrm{i}\pi\nu/2}\,K_{\nu}(-\mathrm{i}z), \qquad
H_{\nu}^{(2)}(z) = -\frac{2}{\pi\,\mathrm{i}}\,
  e^{+\mathrm{i}\pi\nu/2}\,K_{\nu}(+\mathrm{i}z),
\qquad \nu \notin \mathbb{Z},
\label{eq:H1H2_relations_intro}
\end{equation}

\noindent
The evaluation of \(K_{\nu}(z)\) required in the present context is carried out using the robust and fully branch-aware algorithm developed in
Part~II~\cite{Zaghloul_Johnson_II}, where all cases arising from analytic continuation in the complex plane are treated in detail. The same implementation is reused here without modification.

For completeness, we recall the continuation formula that applies when the rotated argument leaves the principal strip \((-\pi/2,\pi/2]\):
\begin{equation}
K_{\nu}(w e^{\pm \mathrm{i}\pi})
  = e^{\mp \mathrm{i}\pi\nu}\,K_{\nu}(w)
    \,\mp\, \mathrm{i}\pi\, I_{\nu}(w).
\label{eq:K_continuation}
\end{equation}
Together with the branch-handling strategy described in Part~II, this relation ensures smooth and accurate behavior for all real \(\nu\) and complex \(z\), including regimes of large order and highly oscillatory arguments.

\noindent
The reflection identity~\eqref{eq:Y_reflection} is numerically efficient and highly accurate for non-integer orders, requiring only two evaluations of $J_\nu(z)$.  The fallback path~\eqref{eq:Y_K_mapping} avoids the $\sin(\pi\nu)$ denominator near integer $\nu$ and remains well-behaved throughout the complex plane. Both paths employ phase stabilization via the integer–fractional split of $\nu$, which ensures that all trigonometric factors are computed with minimal loss of precision even for large $|\nu|$. In practice, the threshold $\tau$ can be tuned empirically to maintain targeted accuracy;  \(\tau \approx \sqrt{\varepsilon}\) has been found reliable for both double and quadruple precision, where \(\varepsilon\) denotes the machine epsilon of the working precision.

\medskip
\noindent
The corresponding \texttt{Fortran} routine \texttt{y\_nu\_of\_z} invokes the
existing routines \texttt{j\_nu\_of\_z} (described above) and
\texttt{k\_nu\_of\_z} from Part~II~\cite{Zaghloul_Johnson_II}. The reflection
identity is used for noninteger orders, while near integer \(\nu\) the code
automatically switches to a branch-corrected mapping based on modified Bessel
functions.

\noindent
The hybrid algorithm is summarized in the following pseudocode:

\begin{algorithm}[H]
\caption{Evaluation of $\Ynu(z)$ for real $\nu$ and complex $z$ (robust reflection with $K$-mapping fallback, using the routines
\texttt{j\_nu\_of\_z}, (described in Section~\ref{sec:J_nu}), \texttt{i\_nu\_of\_z}, and \texttt{k\_nu\_of\_z} from Parts~I and II of this series, respectively)}
\label{alg:Ynu}
\begin{algorithmic}[1]
  \STATE \textbf{Input:} real order $\nu$, complex argument $z_{\mathrm{in}}$.
  \STATE \textbf{Output:} $\Ynu(z_{\mathrm{in}})$ in \texttt{Y}, status \texttt{ierr}.

  \STATE Set \texttt{ierr} $\gets 0$, \texttt{conj\_input} $\gets$ \texttt{false}, $z \gets z_{\mathrm{in}}$.
\STATE \textbf{if} $\Im(z_{\mathrm{in}}) < 0$ \textbf{then} 
       \texttt{conj\_input} $\gets$ \texttt{true}; 
       $z \gets \overline{z_{\mathrm{in}}}$.

  \STATE Compute $n \gets \operatorname{NINT}(\nu)$, $f \gets \nu - n$.
  \STATE If $f > 1/2$ then set $n \gets n+1$, $f \gets f-1$; if $f < -1/2$ then set $n \gets n-1$, $f \gets f+1$.
  \STATE Set $i_{\mathrm{par}} \gets n \bmod 2$.
\STATE \textbf{if} $i_{\mathrm{par}} = 0$ \textbf{then} 
       $\cos(\pi\nu) \gets \cos(\pi f)$, $\sin(\pi\nu) \gets \sin(\pi f)$;
       \textbf{else} 
       $\cos(\pi\nu) \gets -\cos(\pi f)$, $\sin(\pi\nu) \gets -\sin(\pi f)$.

  \STATE Choose threshold $\tau$ (e.g., $\tau \gets \sqrt{\varepsilon}$).

  \STATE \textbf{(A) Reflection identity branch}
  \IF{$|\sin(\pi\nu)| \ge \tau$}
    \STATE Compute $\Jnu(z)$ via \texttt{j\_nu\_of\_z}$(\nu,z)$; if it fails, set \texttt{ierr} and return.
    \STATE Compute $J_{-\nu}(z)$ via \texttt{j\_nu\_of\_z}$(-\nu,z)$; if it fails, set \texttt{ierr} and return.

\STATE Set $\texttt{Y} \gets (\cos(\pi\nu)\Jnu(z) - J_{-\nu}(z)) / \sin(\pi\nu)$.

    \STATE \textbf{goto} finalization.
  \ENDIF

  \STATE \textbf{(B) Fallback: $K$-mapping via $H_{\nu}^{(1)}$}
  \STATE Set $w \gets -\mathrm{i}\,z$.
  \STATE Compute $m_4 \gets n \bmod 4$ and set
  \[
    \texttt{base} \gets
    \begin{cases}
      1, & m_4 = 0,\\
      \mathrm{i}, & m_4 = 1,\\
      -1, & m_4 = 2,\\
      -\mathrm{i}, & m_4 = 3;
    \end{cases}
  \]
  \STATE Set $\phi_{\mathrm{half}} \gets \tfrac{\pi}{2} f$ and
         $\texttt{rot\_half} \gets \cos(\phi_{\mathrm{half}})
           + \mathrm{i}\sin(\phi_{\mathrm{half}})$.
  \STATE Set $e_{-\mathrm{i}\pi\nu/2} \gets \overline{\texttt{base} \cdot \texttt{rot\_half}}$.

   \STATE Set \texttt{need\_conj} $\gets$ true if $\Im(w) < 0$, otherwise false.
  \IF{\texttt{need\_conj}}
    \STATE Compute $K_0 \gets \texttt{k\_nu\_of\_z}(|\nu|, \overline{w})$; if it fails, set \texttt{ierr} and return.
    \STATE Set $K_0 \gets \overline{K_0}$.
  \ELSE
    \STATE Compute $K_0 \gets \texttt{k\_nu\_of\_z}(|\nu|, w)$; if it fails, set \texttt{ierr} and return.
  \ENDIF

   \STATE Form $H^{(1)}_{\nu}(z)\gets -(2\mathrm{i}/\pi)\,e^{-\mathrm{i}\pi\nu/2}\,K_0$.

  \STATE Compute $\Jnu(z)$ via \texttt{j\_nu\_of\_z}$(\nu,z)$; if it fails, set \texttt{ierr} and return.
  \STATE Set $\texttt{Y} \gets -\mathrm{i}\,\bigl(H^{(1)}_{\nu}(z) - \Jnu(z)\bigr)$.

  \STATE \textbf{finalization:}
\STATE \textbf{if} \texttt{conj\_input} \textbf{then} $\texttt{Y} \gets \overline{\texttt{Y}}$.
  \STATE Return $\Ynu(z_{\mathrm{in}}) \gets \texttt{Y}$ and \texttt{ierr}.
\end{algorithmic}
\end{algorithm}


\section{Implementation and Accuracy Verification}
\label{sec:implementation}

The algorithms developed in Section~\ref{sec:algorithm} have been implemented as \texttt{Fortran} subroutines within a unified module supporting both double- and 
quadruple-precision arithmetic. As in earlier parts of this series, the working precision is controlled by the integer parameter \texttt{rk}, which selects the 
appropriate real kind for all computations; this parameter is defined in the auxiliary module \texttt{set\_rk}.

To assess numerical accuracy and benchmark computational performance, we compare the present implementations for the Bessel functions \( \Jnu(z) \) and \(\Ynu(z)\) against two established references: the widely used double-precision implementation of Algorithm~644~\cite{Amos1986}, and the more recent Algorithm~912~\cite{Kodama2011}, which supports both double- and quadruple-precision arithmetic as well as complex orders.

\noindent
It should be emphasized that Algorithm~912 does not provide uniform coverage of the \((\nu,z)\)-plane. In particular, its applicability is restricted for limited ranges of values of \(|\nu|\) and \(|z|\), as well as in certain regions of the complex plane. This limitation is directly reflected in the numerical comparisons presented below, where Algorithm~912 fails to produce valid outputs over wide regions of the domain representable in double-precision arithmetic, while the present algorithm continues to deliver stable and accurate results.

High-precision reference values were generated using \texttt{Maple}'s arbitrary-precision special-function library and independently validated with the \texttt{mpmath} variable-precision package~\cite{mpmath}. Relative errors were evaluated over extensive grids of real orders and complex arguments \((\nu, z)\), spanning both small and large parameter regimes.

For each test point, the real order \(\nu\), the real part \(x=\Re z\), and the
imaginary part \(y=\Im z\) were sampled primarily on logarithmic scales. The reported relative errors are defined
componentwise by
\begin{equation}
\label{eq:rel_error}
E_{\Re}=\frac{|\Re f_{\rm comp}-\Re f_{\rm ref}|}{|\Re f_{\rm ref}|},
\qquad
E_{\Im}=\frac{|\Im f_{\rm comp}-\Im f_{\rm ref}|}{|\Im f_{\rm ref}|}.
\end{equation}
When the reference component is zero or below the safe relative-error threshold, that point is excluded from the componentwise relative-error summary. This avoids artificially large relative errors caused solely by division by a vanishing reference component.

\subsection{Accuracy Verification in Double Precision} \label{sec:acc_dp}
In this subsection, we assess the accuracy of the present algorithms in double precision. Since Algorithm~644 is restricted to nonnegative real orders and double-precision arithmetic, comparisons with this method are limited to positive orders. Algorithm~912, by contrast, supports complex orders and complex arguments; therefore, the comparisons with Algorithm~912 include both positive and negative real orders, while also illustrating its restricted domain of applicability.
 
\subsubsection{Bessel Functions of the First Kind, \( \Jnu(z) \). } \label{sec:acc_dp_j}
To evaluate the accuracy of the double-precision implementation of \(\Jnu(z)\), we used a grid of 281{,}826 test points distributed primarily on a logarithmic scale across the \(\nu\), \(x\) and \(y\) domains. 
All reference values at these points were verified to lie within the representable range of IEEE double precision, ensuring that the observed errors reflect algorithmic behavior rather than floating-point limitations.

\begin{figure}[htbp]
    \centering
    \subfloat[]{\includegraphics[width=0.48\linewidth]{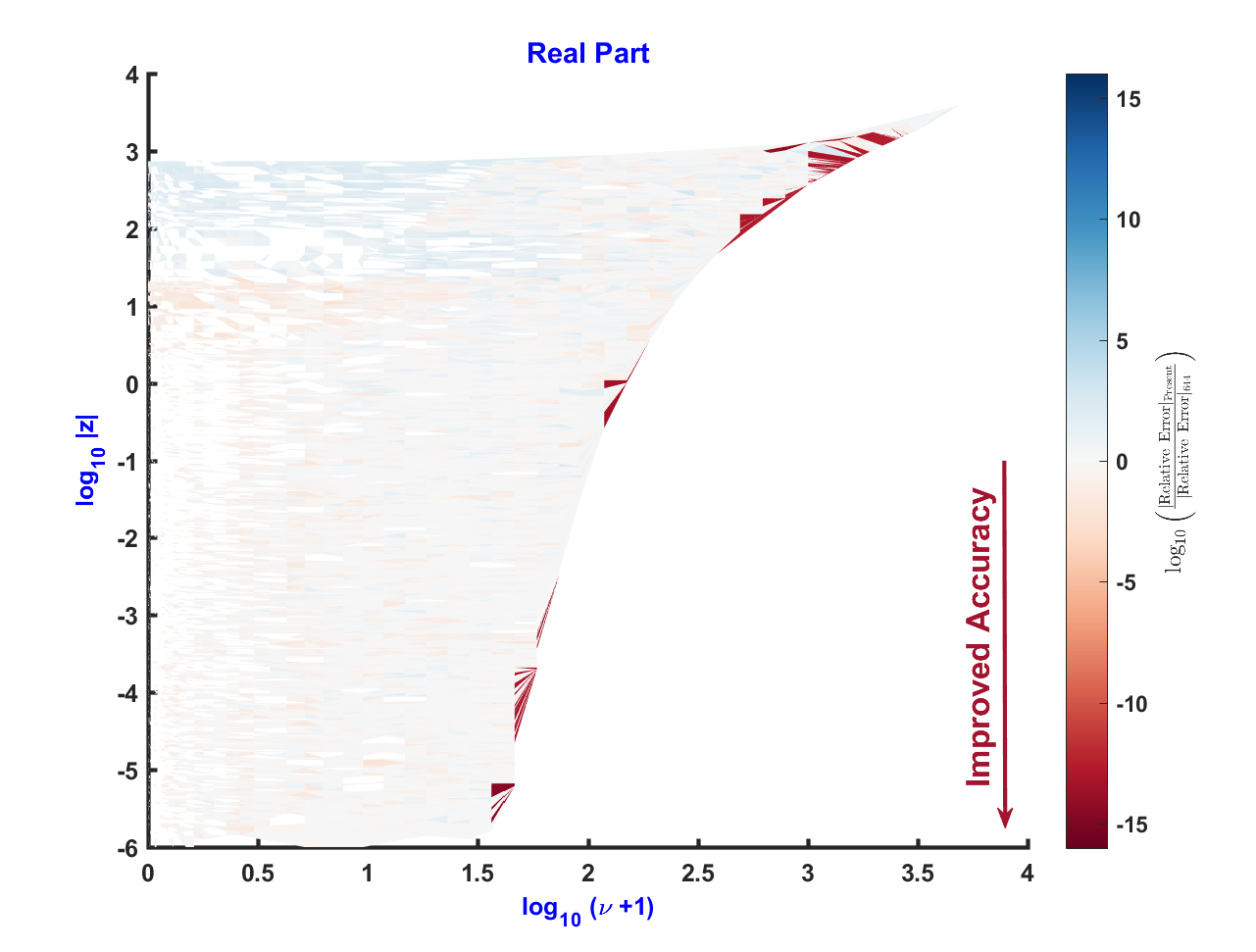}}
    \subfloat[]{\includegraphics[width=0.48\linewidth]{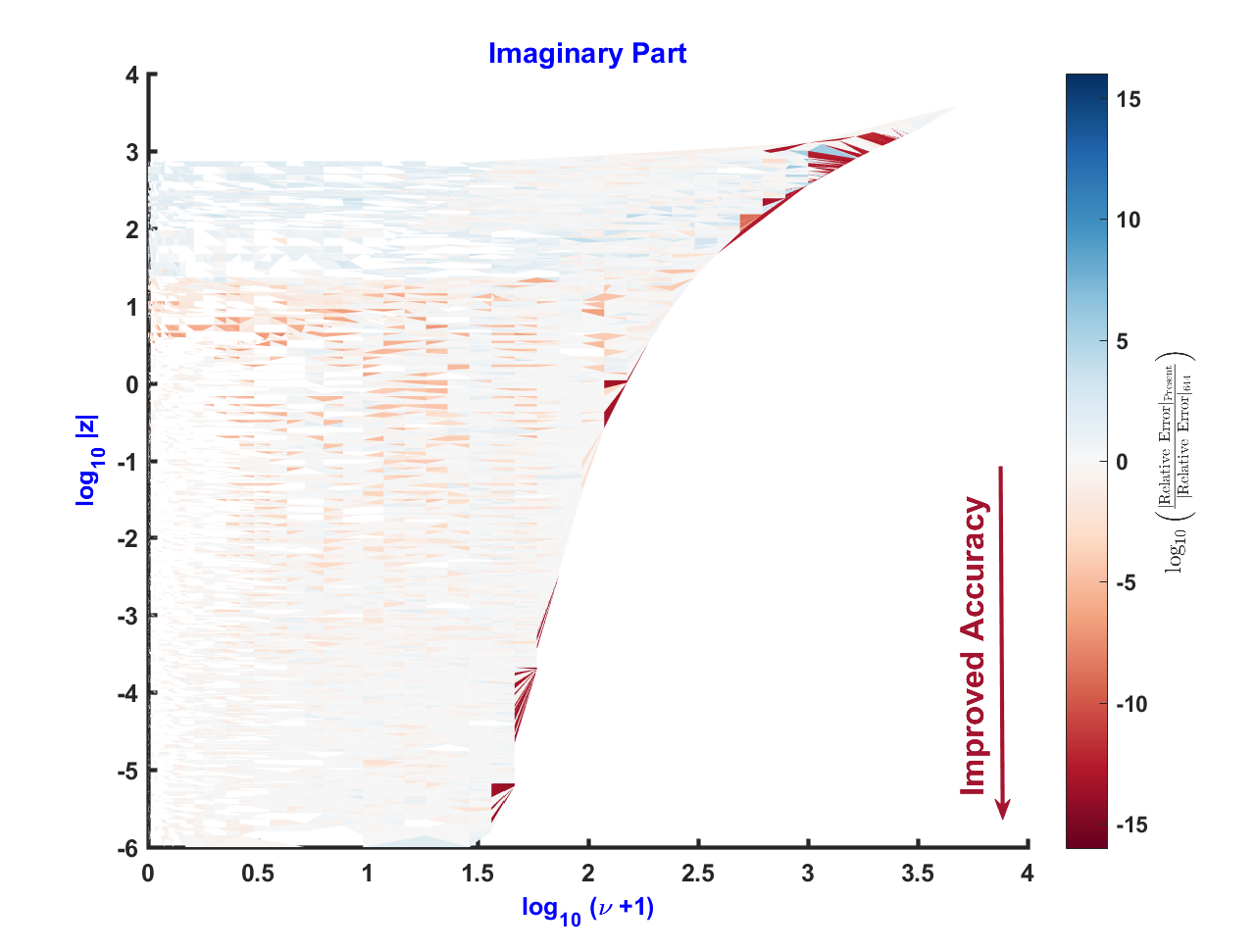}}
    \caption{Colormap plots of the base-10 logarithm of the \textit{ratio} of relative errors (present/Algorithm~644) in computing the real part~(a) and imaginary part~(b) of \(\Jnu(z)\) for the 281{,}826 test points, using \texttt{Maple} values as the reference.}
    \label{fig:acc_rel_J}
\end{figure}

Figure~\ref{fig:acc_rel_J} presents colormap plots of the base-10 logarithm of the ratio of the relative error of the present implementation to that of Algorithm~644, for the real part~(a) and imaginary part~(b) of \(\Jnu(z)\), with \texttt{Maple} serving as a common reference. As in earlier parts of this series, Algorithm~644 fails over significant regions of the domain due to overly conservative overflow and underflow thresholds. In such cases, the routine returns \((0,0)\), corresponding to a \(100\%\) relative error; these failures appear as red regions in Fig.~\ref{fig:acc_rel_J}. By contrast, the present algorithm continues to provide accurate results throughout these regions. In the regions where Algorithm~644 returns computed values, its accuracy is comparable to that of the present method.

\begin{figure}[htbp]
    \centering
    \subfloat[]{\includegraphics[width=0.48\linewidth]{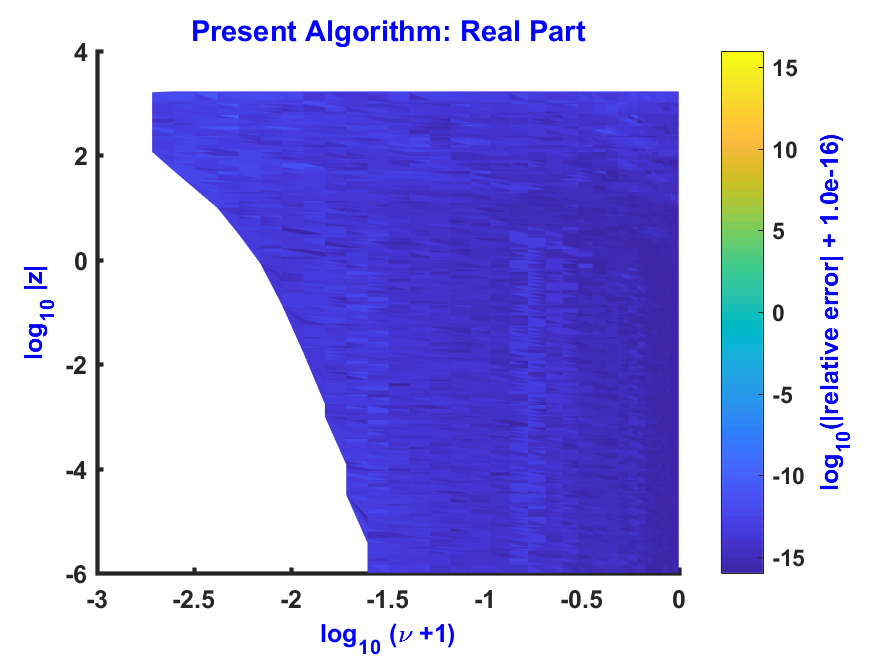}}
    \subfloat[]{\includegraphics[width=0.48\linewidth]{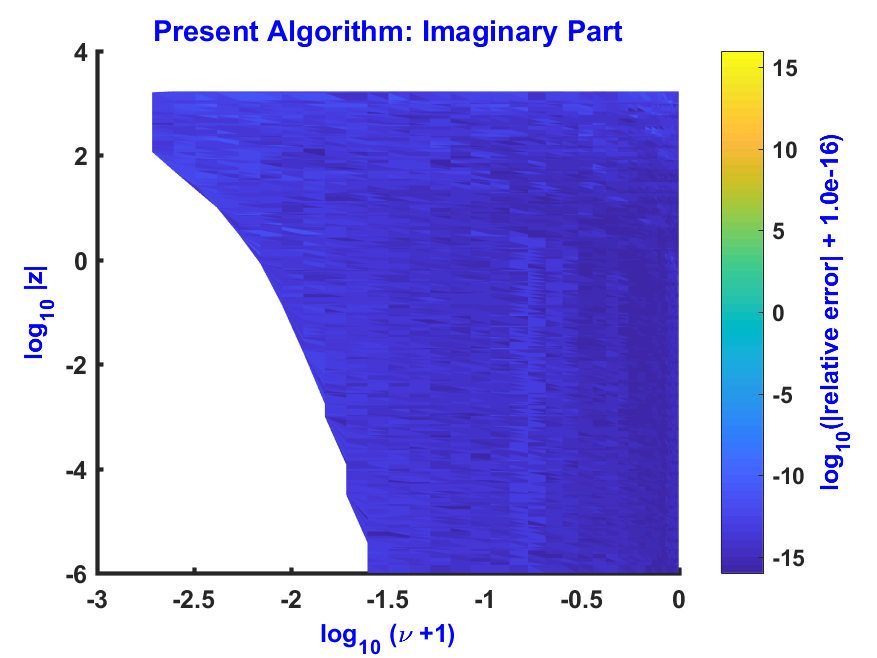}}
    \caption{Colormap plots of the base-10 logarithm of the relative error $(\,|\text{present} - \text{reference}| \,/\, |\text{reference}|\,)$
 in computing the real part~(a) and imaginary part~(b)
    of \(J_{-|\nu|}(z)\) for the double-precision test grid. \texttt{Maple} values are
    used as the reference.}
    \label{fig:acc_j_dbl_neg_nu}
\end{figure}

As noted above, for negative orders Algorithm~644 immediately returns an error flag and does not perform any computation. Consequently, accuracy in this regime is assessed solely against high-precision \texttt{Maple} reference values. As shown in Fig.~\ref{fig:acc_j_dbl_neg_nu}, the present algorithm computes \(J_{-|\nu|}(z)\) accurately across the tested domain, in agreement with the theoretical symmetry \(J_{-|n|}(z)=(-1)^{|n|}J_{|n|}(z)\) for integer orders and with its analytic continuation to noninteger negative orders.

\begin{figure}[htbp]
    \centering
    \subfloat[]{\includegraphics[width=0.48\linewidth]{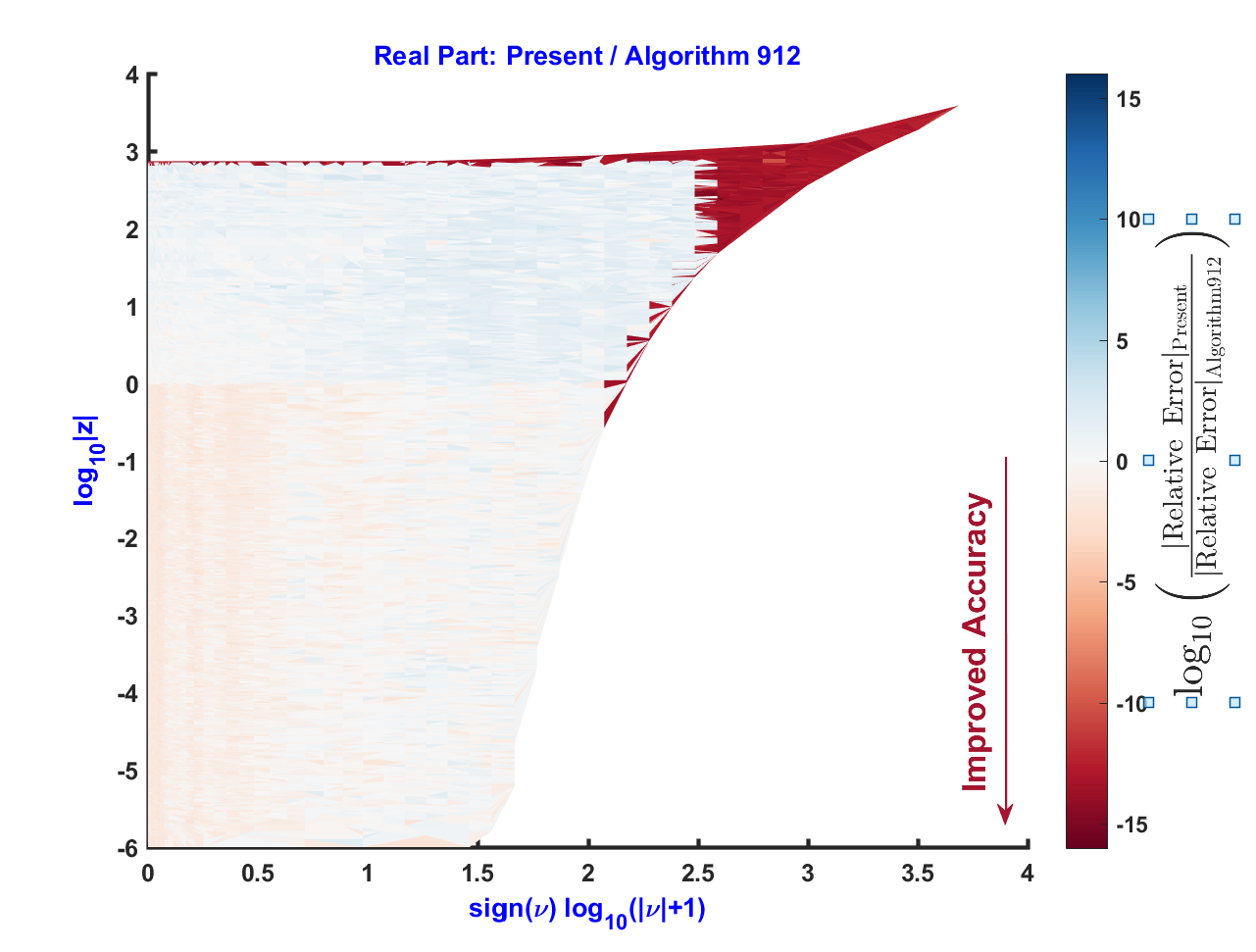}}
    \subfloat[]{\includegraphics[width=0.48\linewidth]{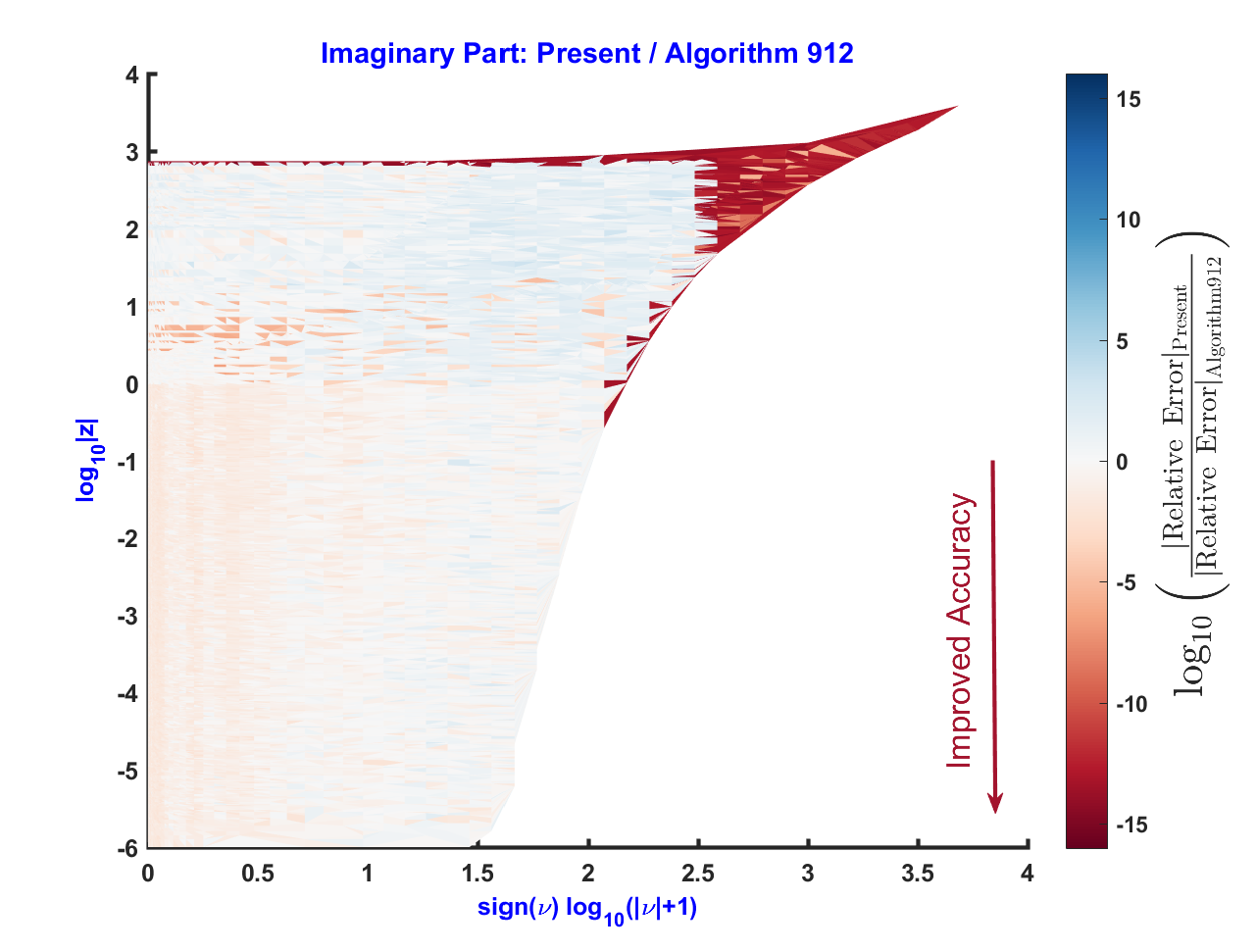}}
    \caption{Colormap plots of the base-10 logarithm of the \textit{ratio} of relative errors (present/Algorithm~912) in computing the real part~(a) and imaginary part~(b) of \(\Jnu(z)\) for the 281{,}826 test points, using \texttt{Maple} values as the reference.}
    \label{fig:acc_rel_J_912}
\end{figure}

\begin{figure}[htbp]
    \centering
    \subfloat[]{\includegraphics[width=0.48\linewidth]{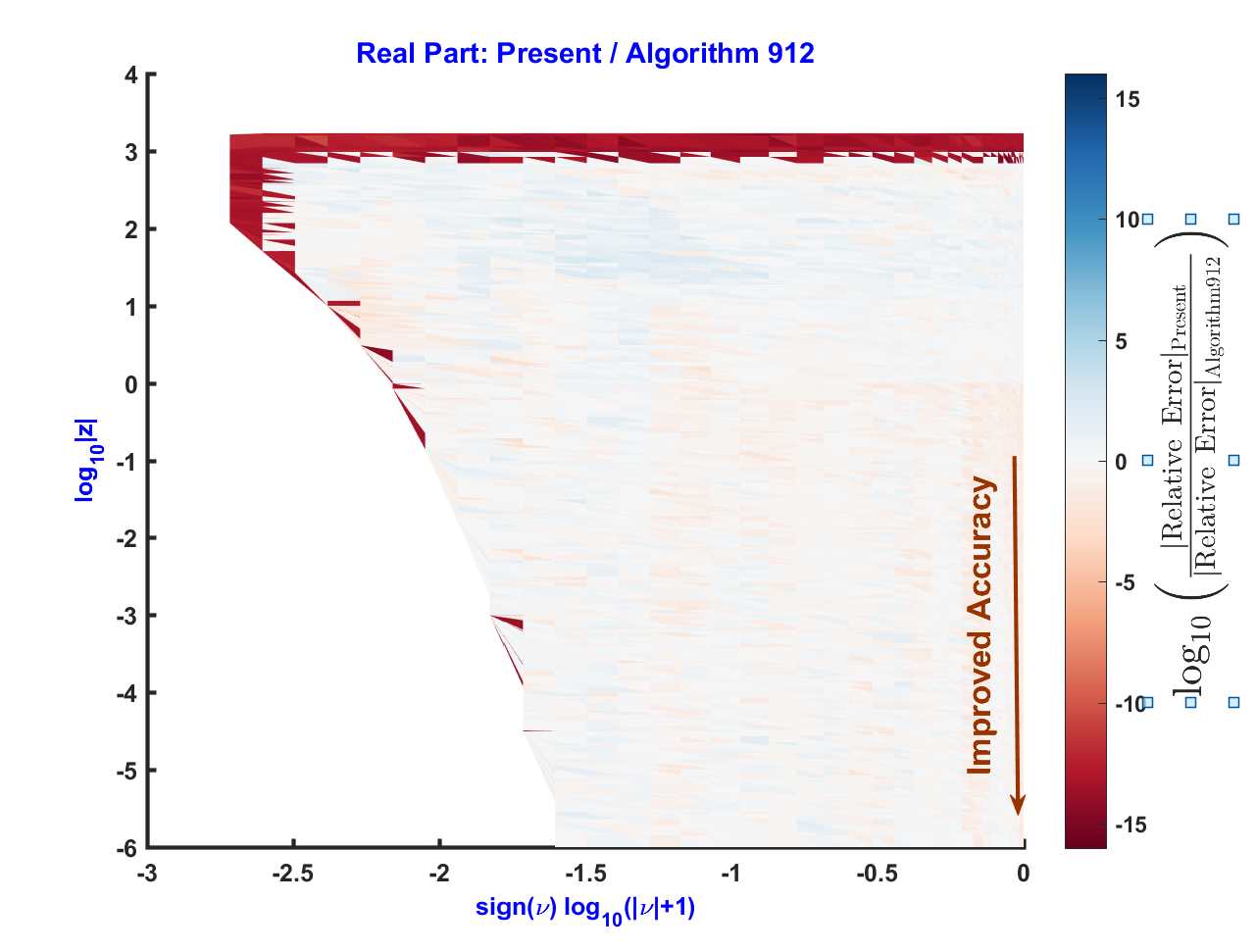}}
    \subfloat[]{\includegraphics[width=0.48\linewidth]{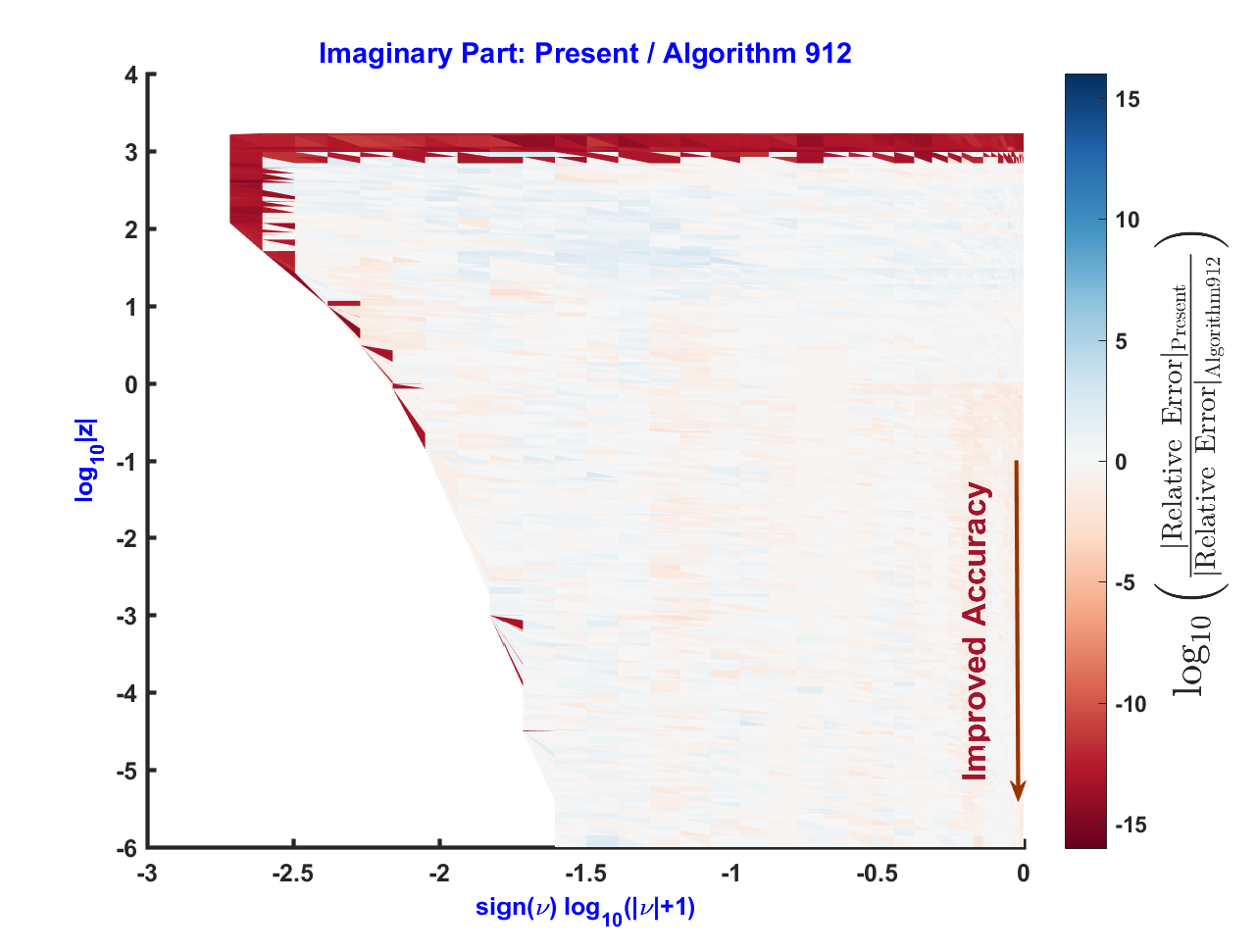}}
    \caption{Colormap plots of the base-10 logarithm of the \textit{ratio} of relative errors (present/Algorithm~912) in computing the real part~(a) and imaginary part~(b) of\(J_{-|\nu|}(z)\) for the 281{,}826 test points, using \texttt{Maple} values as the reference.}
    \label{fig:acc_rel_J_912_neg}
\end{figure}

Figures~\ref{fig:acc_rel_J_912} and~\ref{fig:acc_rel_J_912_neg} compare the accuracy of the present implementation with Algorithm~912, using high-precision \texttt{Maple} results as the reference. The comparison confirms that Algorithm~912 is restricted to a limited region of applicability, as documented in its original description. Within this admissible domain, both methods exhibit comparable levels of accuracy. Outside this region, however, Algorithm~912 is no longer applicable, whereas the present algorithm continues to provide stable and accurate results over a substantially broader domain.

\subsubsection{Bessel Functions of the Second Kind, \( \Ynu(z) \).} \label{sec:acc_dp_Y}

The accuracy of the double-precision implementation of \(\Ynu(z)\) was examined using a grid of 281{,}935 test points, likewise distributed primarily on a
logarithmic scale across the \(\nu\), \(x\), and \(y\) domains. As before, all reference values were verified to lie within the representable range of IEEE
double-precision arithmetic, ensuring that the resulting error measurements reflect the numerical behavior of the algorithms rather than floating-point limitations.

\begin{figure}[htbp]
    \centering
    \subfloat[]{\includegraphics[width=0.48\linewidth]{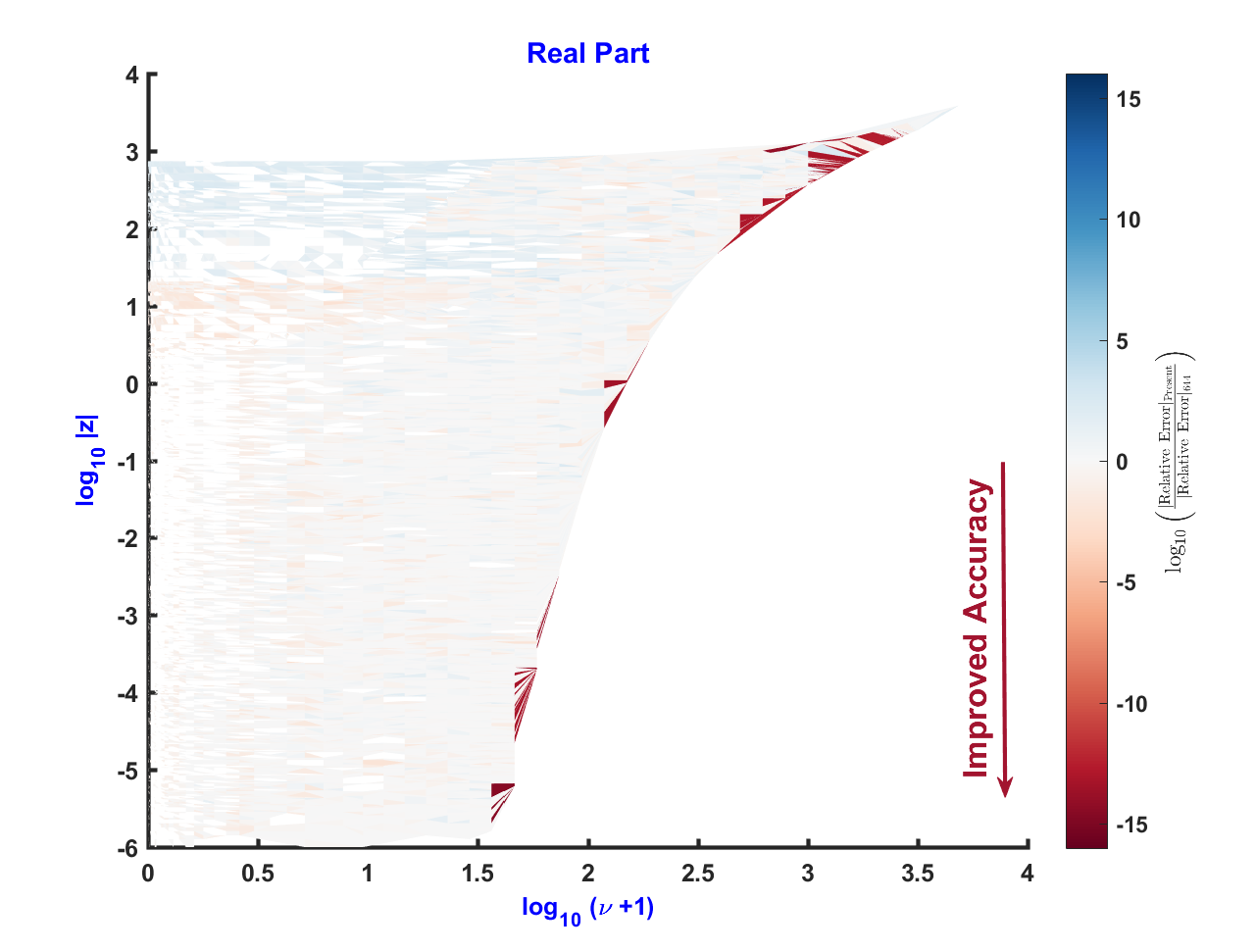}}
    \subfloat[]{\includegraphics[width=0.48\linewidth]{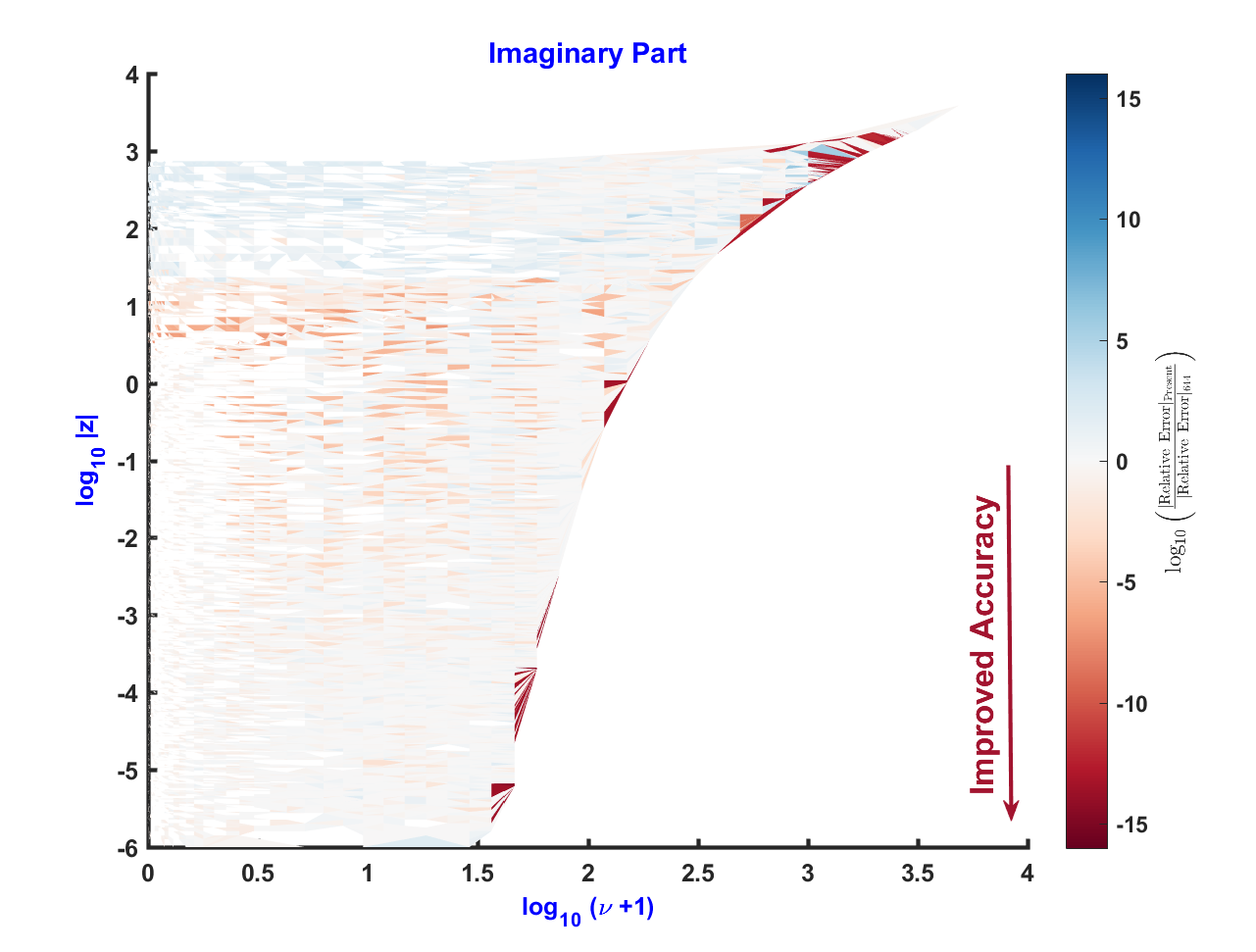}}
    \caption{Colormap plots of the base-10 logarithm of the \textit{ratio} of relative errors (present/Algorithm~644) in computing the real part~(a) and imaginary part~(b) of \(\Ynu(z)\) in double precision. \texttt{Maple} values are used as the reference.}
    \label{fig:acc_rel_Y}
\end{figure}

\noindent
Similar to the case of $\Jnu(z)$, Fig.~\ref{fig:acc_rel_Y} shows the base-10 logarithm of the \textit{ratio} of relative errors between the present implementation and Algorithm~644 for the real part~(a) and imaginary part~(b) of \(\Ynu(z)\), using \texttt{Maple} as the common reference.

The behavior of Algorithm~644 for \(\Ynu(z)\) reflects the same structural limitations observed for \(\Jnu(z)\). The routine terminates prematurely across large regions of the \((\nu,z)\)-plane due to conservative overflow and underflow criteria. These failures again manifest as red regions in Fig.~\ref{fig:acc_rel_Y}. On the other hand, the present algorithm remains stable over the entire test grid, including parameter regimes where Algorithm~644 provides no usable output. In regions of the domain where Algorithm~644 does return values, its accuracy is broadly consistent with that of the present implementation.

\begin{figure}[htbp]
    \centering
    \subfloat[]{\includegraphics[width=0.48\linewidth]{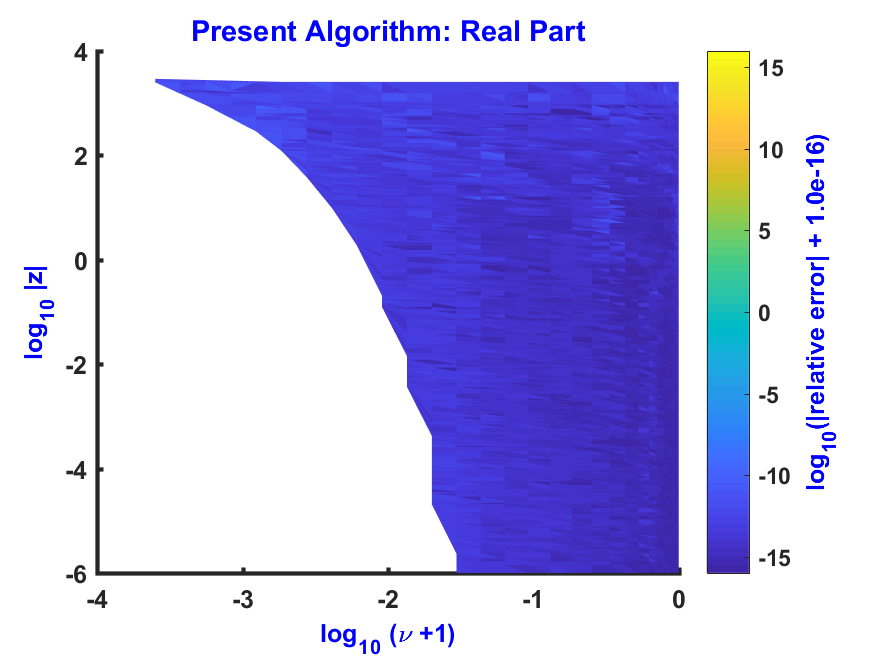}}
    \subfloat[]{\includegraphics[width=0.48\linewidth]{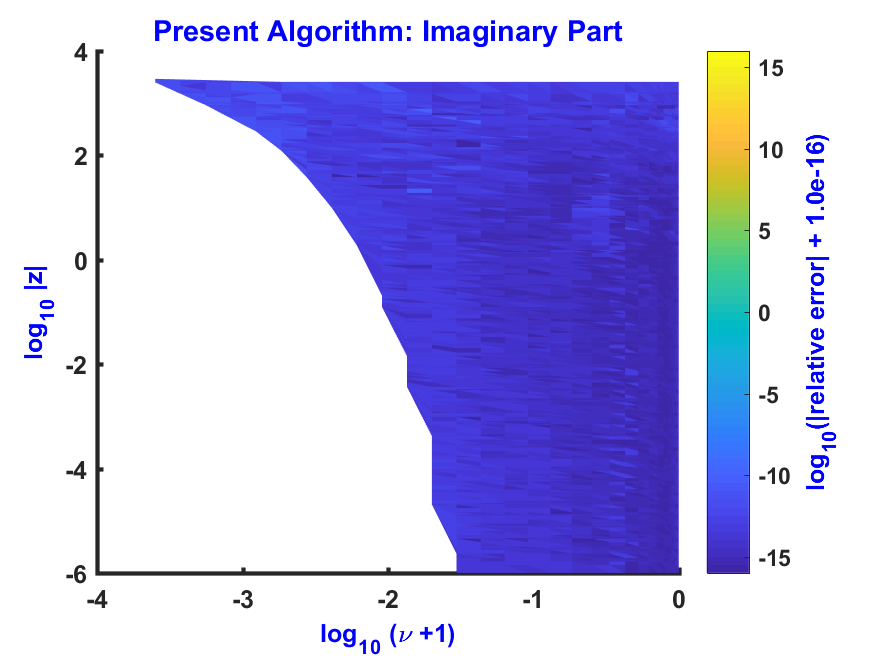}}
    \caption{Colormap plots of the base-10 logarithm of the relative error
   $(\,|\text{present} - \text{reference}| \,/\, |\text{reference}|\,)$ in computing the real part~(a) and imaginary part~(b)
    of \(Y_{-|\nu|}(z)\) in double precision. \texttt{Maple} values are used as the
    reference.}
    \label{fig:acc_rel_Y_dp_neg_nu}
\end{figure}

\begin{figure}[htbp]
    \centering
    \subfloat[]{\includegraphics[width=0.48\linewidth]{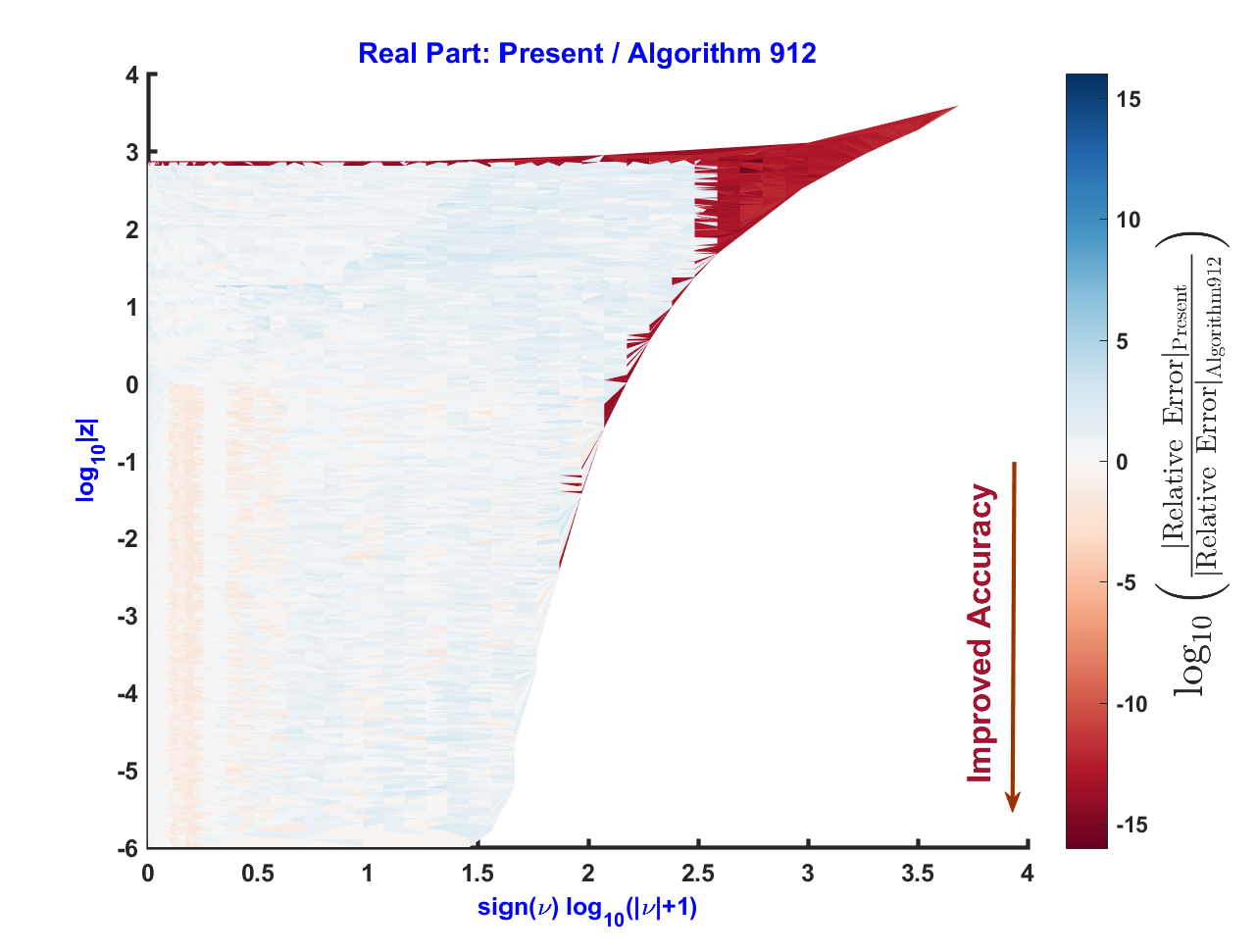}}
    \subfloat[]{\includegraphics[width=0.48\linewidth]{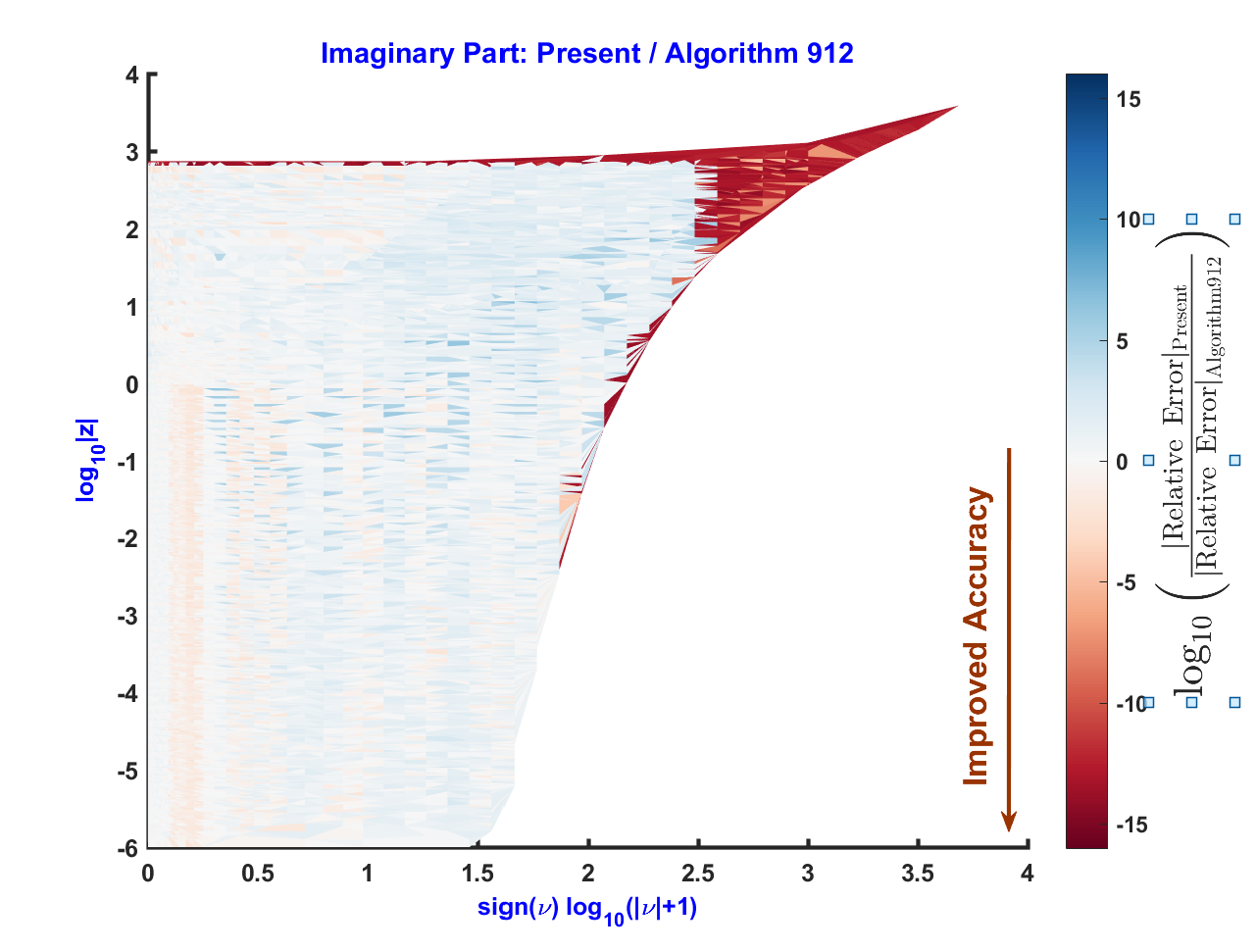}}
    \caption{Colormap plots of the base-10 logarithm of the \textit{ratio} of relative errors (present/Algorithm~912) in computing the real part~(a) and imaginary part~(b) of \(\Ynu(z)\) in double precision. \texttt{Maple} values are used as the reference.}
    \label{fig:acc_rel_Y_912}
\end{figure}

\begin{figure}[htbp]
    \centering
    \subfloat[]{\includegraphics[width=0.48\linewidth]{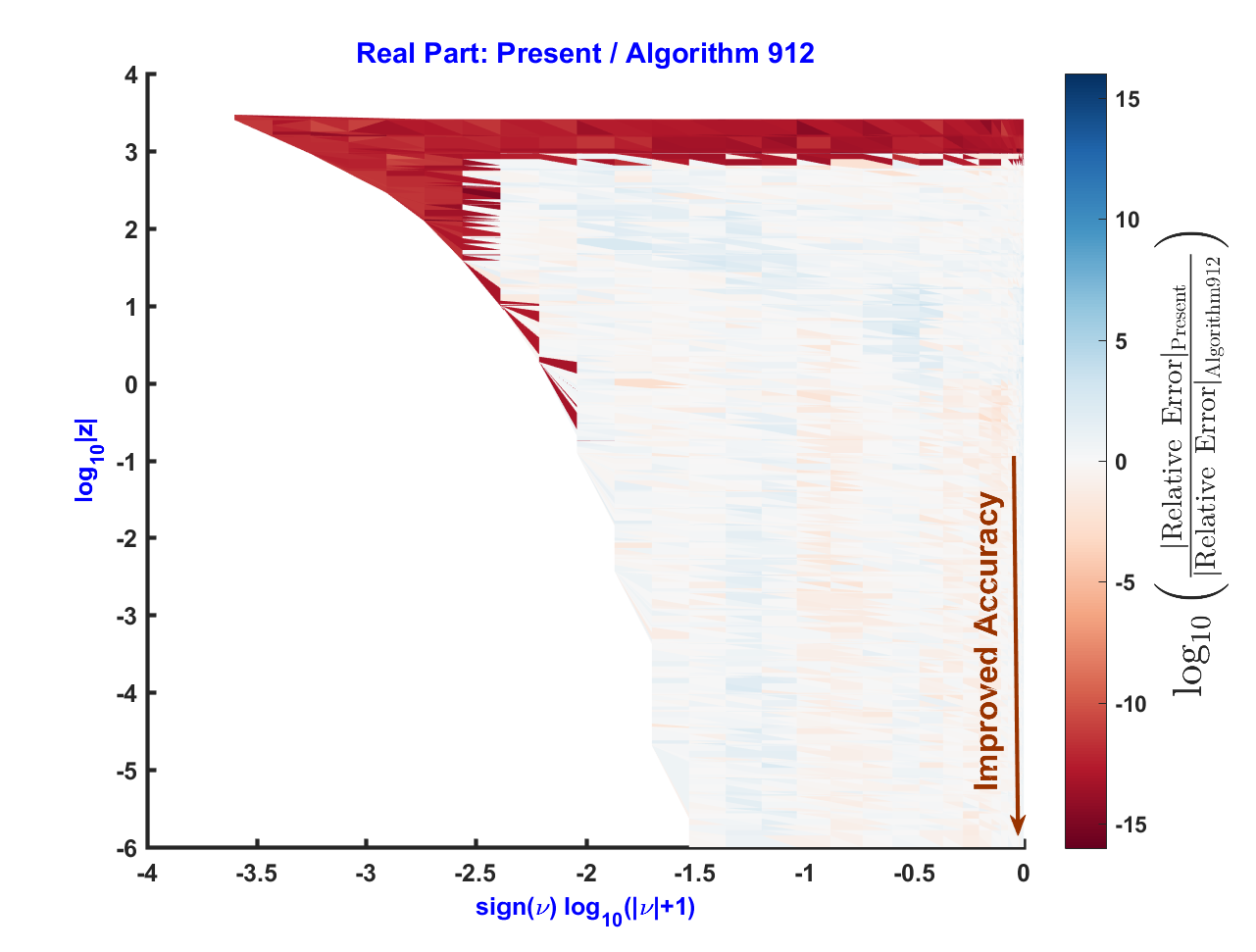}}
    \subfloat[]{\includegraphics[width=0.48\linewidth]{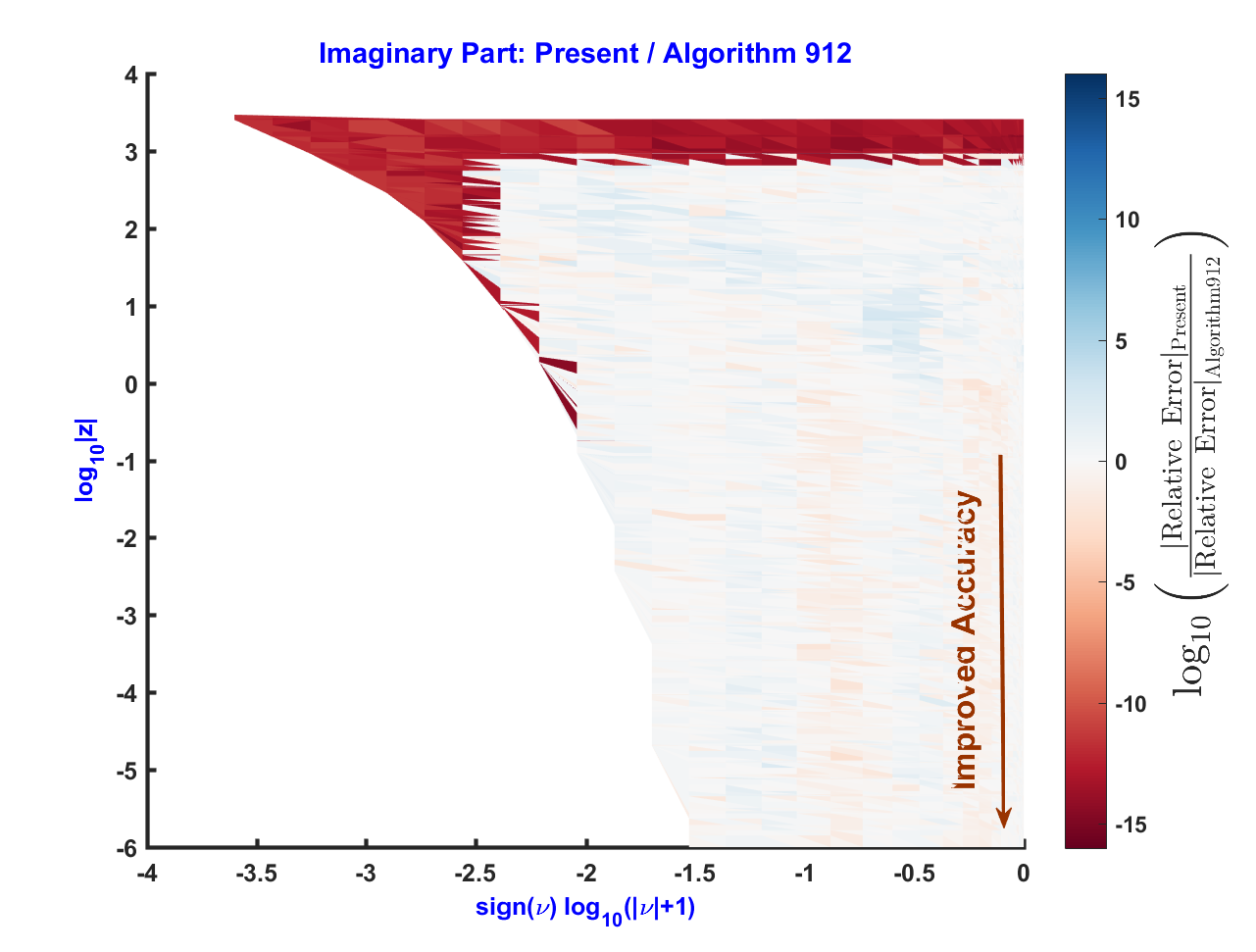}}
    \caption{Colormap plots of the base-10 logarithm of the \textit{ratio} of relative errors (present/Algorithm~912) in computing the real part~(a) and imaginary part~(b) of \(Y_{-|\nu|}(z)\) in double precision, using \texttt{Maple} values as the reference.}
    \label{fig:acc_rel_Y_neg_nu_912}
\end{figure}

\noindent
As with the first-kind Bessel function, Algorithm~644 does not support negative orders of \(\Ynu\), and accuracy in this regime must therefore be assessed 
directly against high-precision \texttt{Maple} reference data. 

Figure~\ref{fig:acc_rel_Y_dp_neg_nu} presents the corresponding relative-error colormaps for \(Y_{-|\nu|}(z)\). The results demonstrate that the present algorithm is accurate across the full tested domain, including the neighborhoods where analytic continuation introduces additional numerical challenges. This confirms the stability and reliability of the proposed methods for both positive and negative orders of \(\Ynu(z)\).

\noindent
A similar trend is observed for the Bessel function of the second kind. As shown in Figures~\ref{fig:acc_rel_Y_912} and~\ref{fig:acc_rel_Y_neg_nu_912}, Algorithm~912 achieves accuracy comparable to the present implementation within its admissible parameter range. Nevertheless, its applicability remains limited, particularly for large orders and arguments, whereas the present algorithm maintains robustness across the entire tested domain.


\subsection{Accuracy Verification in Quadruple Precision} \label{sec:acc_qp}
In this section we assess the performance of the proposed algorithms in native quadruple precision. By increasing the significand from approximately 16 to approximately 34 decimal digits, quadruple precision significantly expands the reliably computable domain in both \(\nu\) and \(z\), typically by about an order of magnitude. This increased precision also mitigates the effects of round-off error, particularly in numerically delicate operations such as analytic continuation and cancellation-prone combinations.

Since Algorithm~644 is restricted to double-precision arithmetic, accuracy verification in quadruple precision is performed primarily against high-precision \texttt{Maple} reference data, and, where applicable, against Algorithm~912~\cite{Kodama2011}, which provides quadruple-precision support over a limited parameter range.

\subsubsection{Bessel Functions of the First Kind, \texorpdfstring{$\Jnu(z)$}{Jnu(z)}} \label{sec:acc_qp_J}
To assess the accuracy of the quadruple-precision implementation of \(\Jnu(z)\), a new test grid of 280{,}636 points was generated over the enlarged \(\nu\)–\(|z|\) domain enabled by the increased dynamic range of quadruple precision. This grid includes values of \(\nu\) and \(z\) that extend well beyond the limits of reliable computation in double precision, thereby providing a more stringent test of numerical stability and accuracy.

\begin{figure}[htbp]
    \centering
    \subfloat[]{\includegraphics[width=0.48\linewidth]{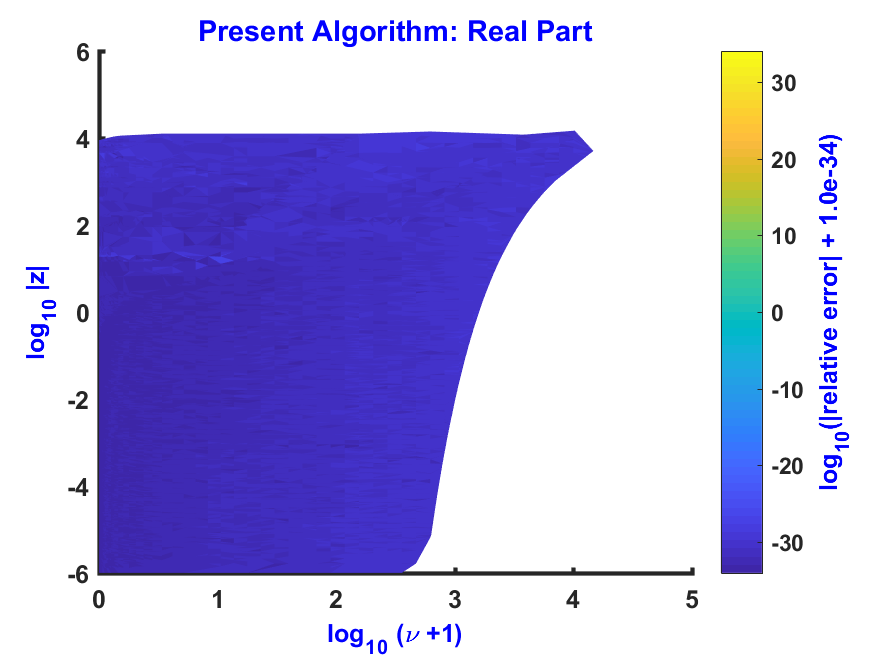}}
    \subfloat[]{\includegraphics[width=0.48\linewidth]{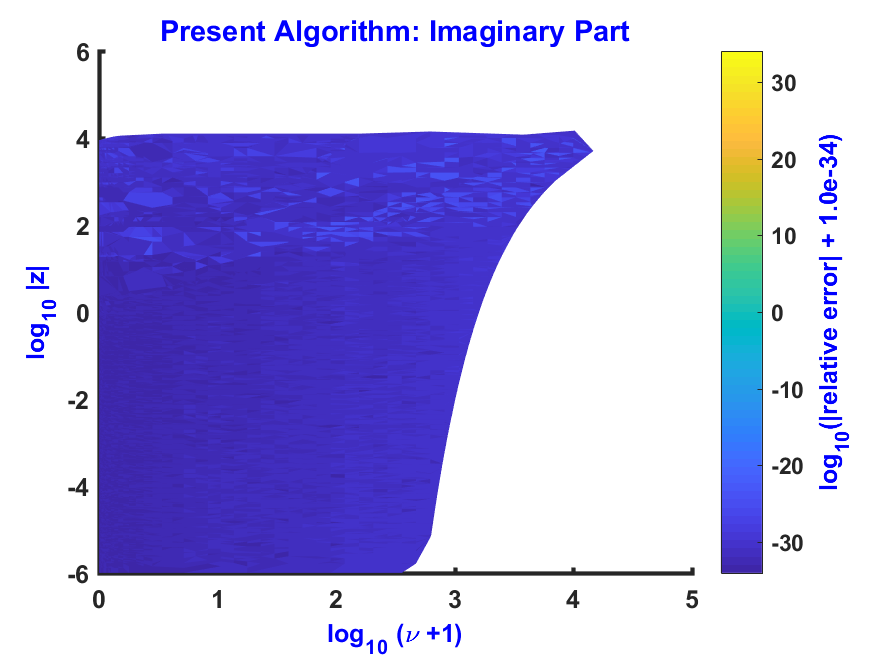}}
    \caption{Colormap plots of the base-10 logarithm of the componentwise relative error in computing the real
    part~(a) and imaginary part~(b) of \(\Jnu(z)\) in quadruple precision
    over the expanded test grid, using \texttt{Maple} values as the reference.}
    \label{fig:qp_err_J_nu}
\end{figure}

Figure~\ref{fig:qp_err_J_nu} presents colormap plots of the relative error in computing the real part~(a) and imaginary part~(b) of \(\Jnu(z)\), using high-precision \texttt{Maple} values as the reference. The results demonstrate that the present algorithm maintains uniformly high accuracy across the entire extended domain, with relative errors reduced by several orders of magnitude compared to those observed in double precision. These results confirm the effectiveness of the proposed hybrid strategy in controlling error propagation even in regimes involving large orders and highly oscillatory complex arguments.

\begin{figure}[htbp]
    \centering
    \subfloat[]{\includegraphics[width=0.48\linewidth]{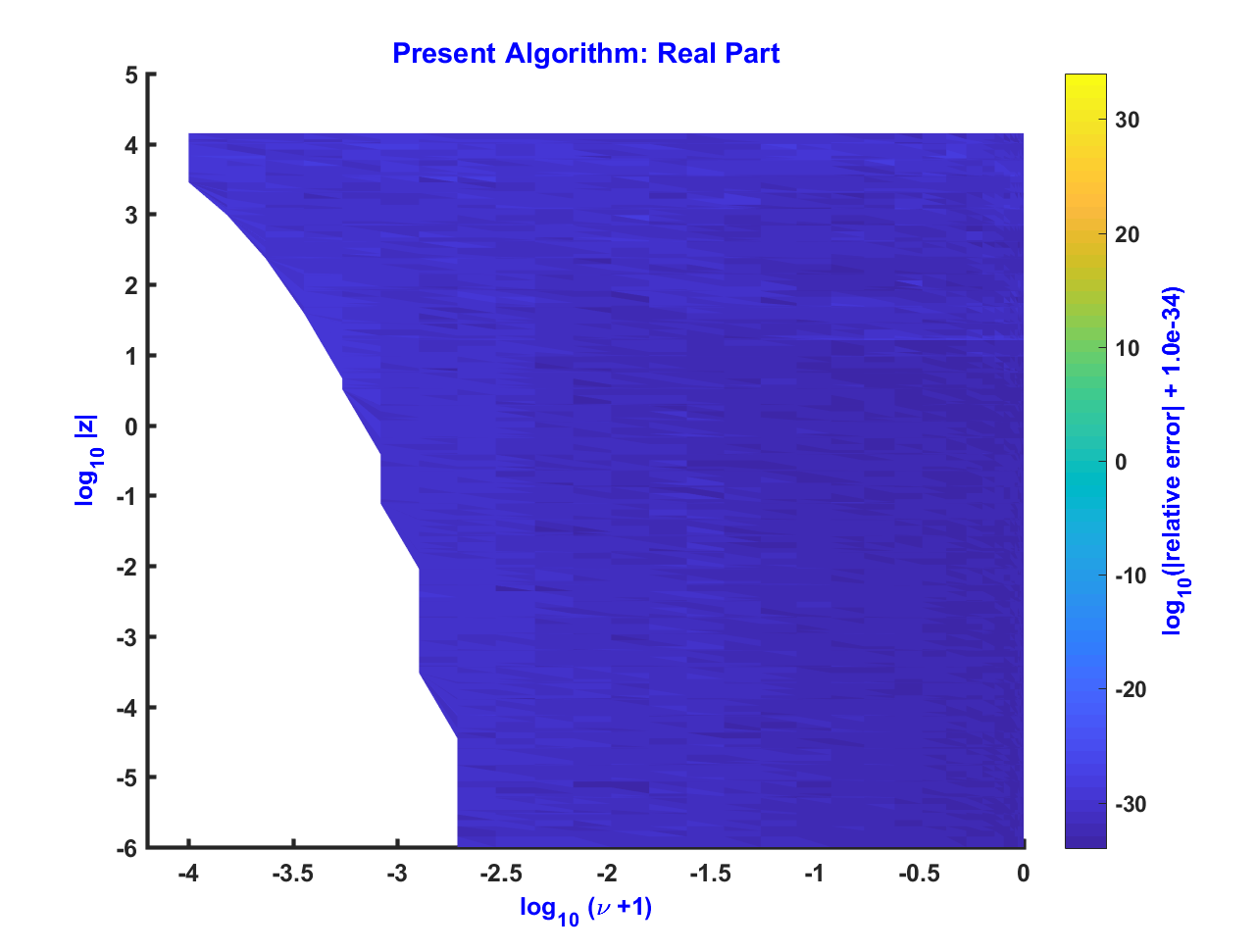}}
    \subfloat[]{\includegraphics[width=0.48\linewidth]{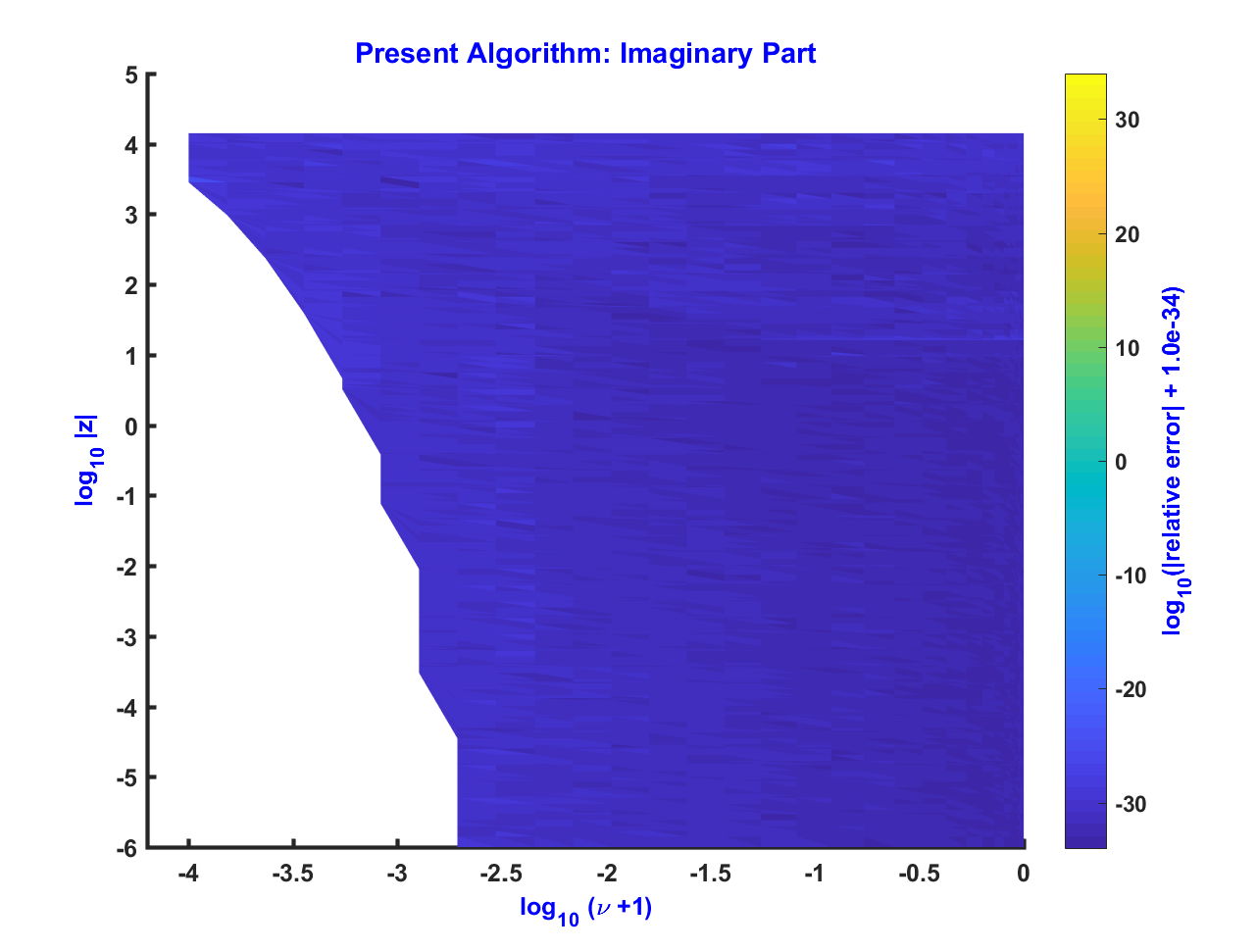}}
    \caption{Colormap plots of the base-10 logarithm of the relative error in computing the real
    part~(a) and imaginary part~(b) of \(J_{-|\nu|}(z)\) in quadruple precision,
    using \texttt{Maple} values as the reference.}
    \label{fig:qp_err_j_neg_nu}
\end{figure}

To further examine the robustness of the method, we consider negative real orders using a separate test grid of 278{,}804 points spanning different sign combinations of the real and imaginary parts of \(z\). As shown in Fig.~\ref{fig:qp_err_j_neg_nu}, the present implementation preserves high accuracy throughout the domain, in excellent agreement with \texttt{Maple} reference values, thereby confirming the stability of the analytic continuation procedures employed.

\begin{figure}[htbp]
    \centering
    \subfloat[]{\includegraphics[width=0.48\linewidth]{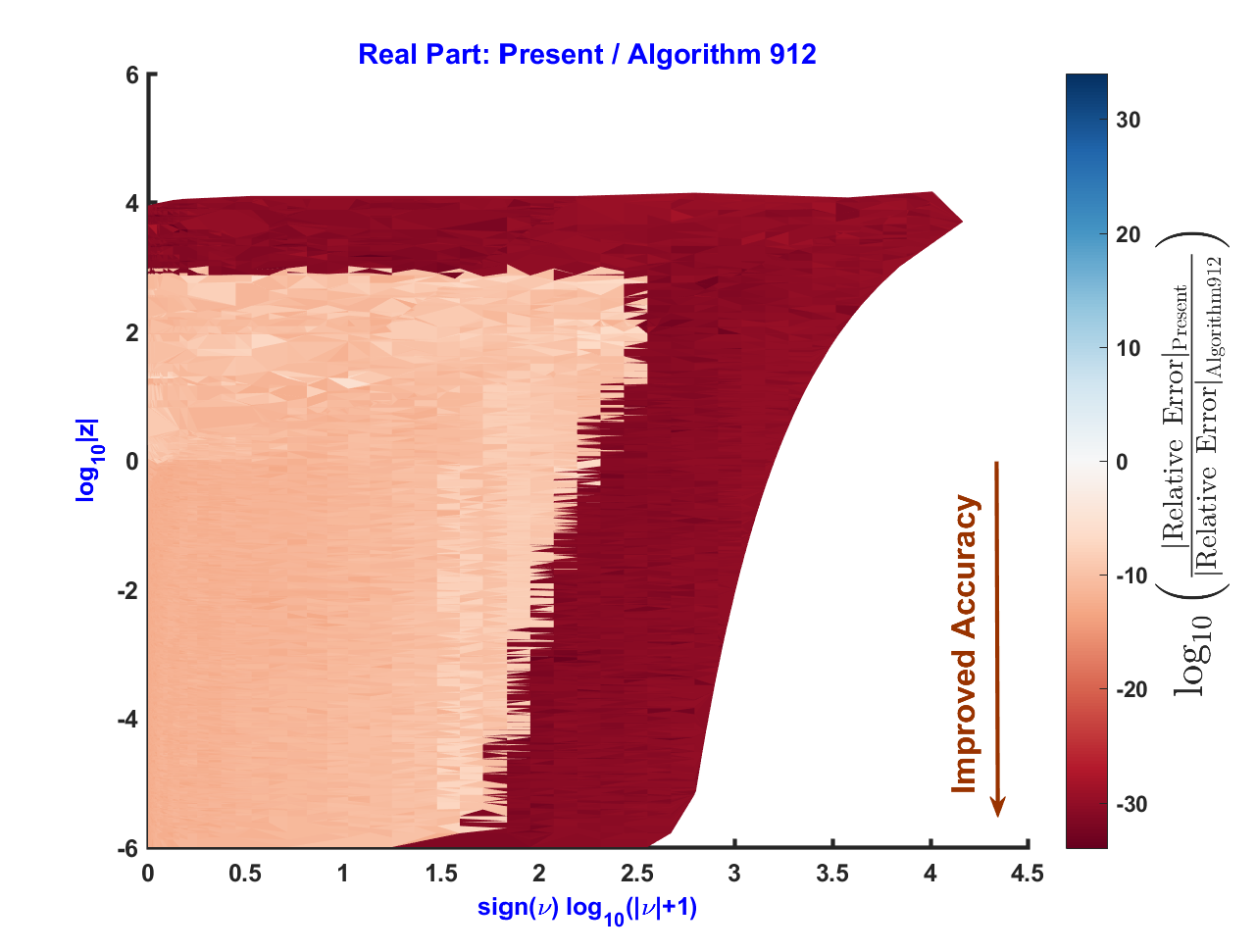}}
    \subfloat[]{\includegraphics[width=0.48\linewidth]{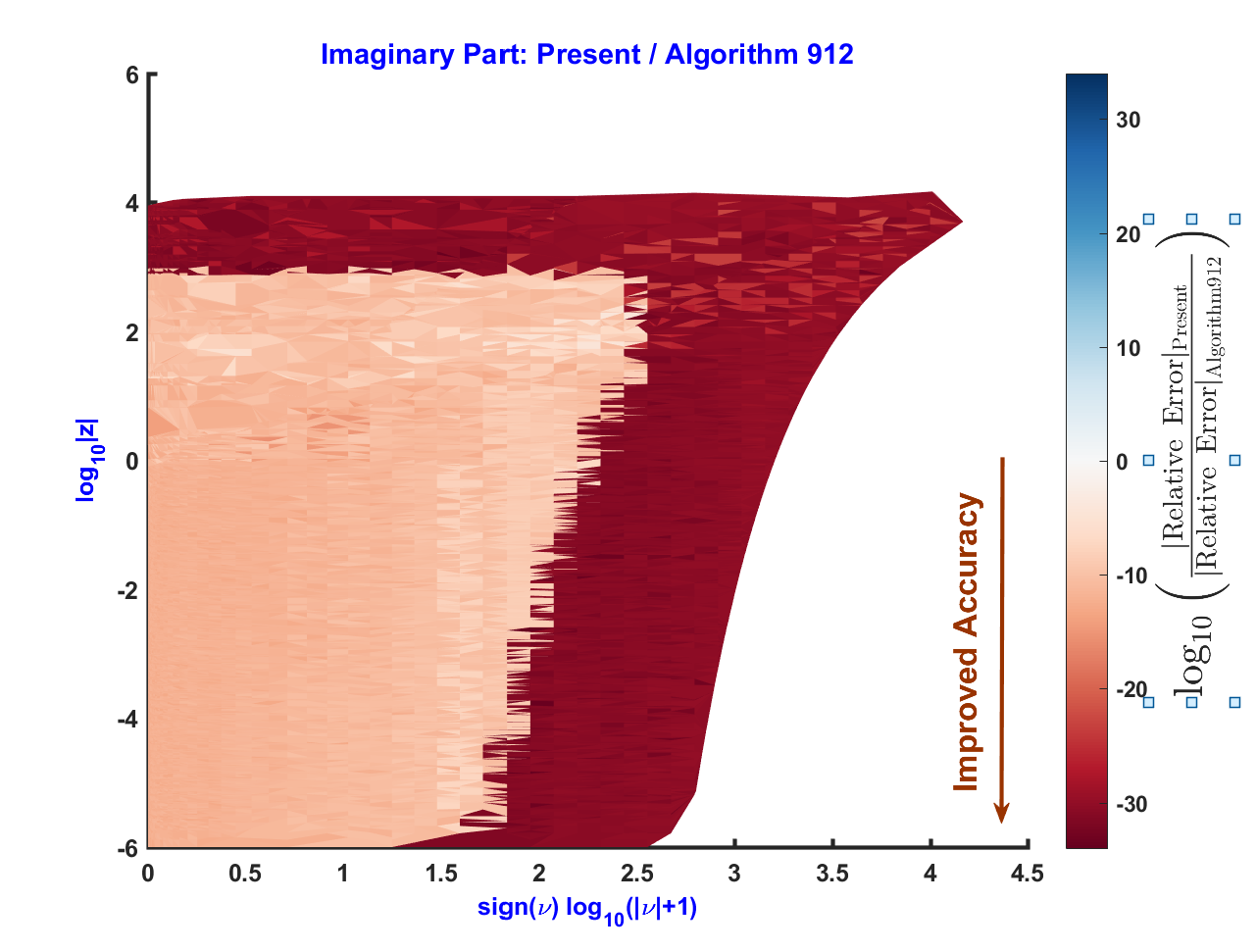}}
    \caption{Colormap plots of the base-10 logarithm of the \textit{ratio} of componentwise
relative errors (present/Algorithm~912) in computing the real part~(a) and
imaginary part~(b) of \(J_{\nu}(z)\) in quadruple precision, using
\texttt{Maple} values as the reference.}
    \label{fig:qp_j_nu_re_912}
\end{figure}

\begin{figure}[htbp]
    \centering
    \subfloat[]{\includegraphics[width=0.48\linewidth]{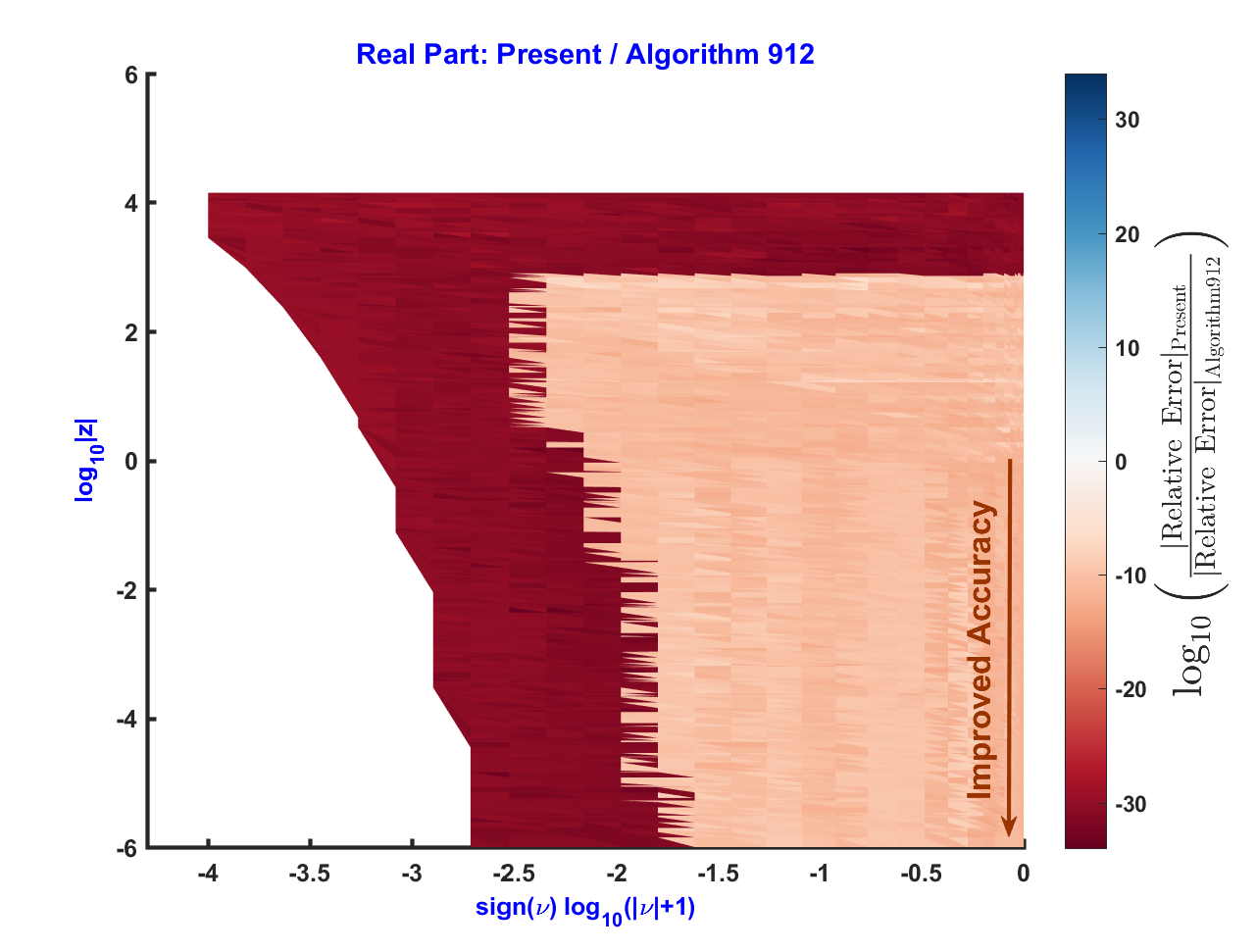}}
    \subfloat[]{\includegraphics[width=0.48\linewidth]{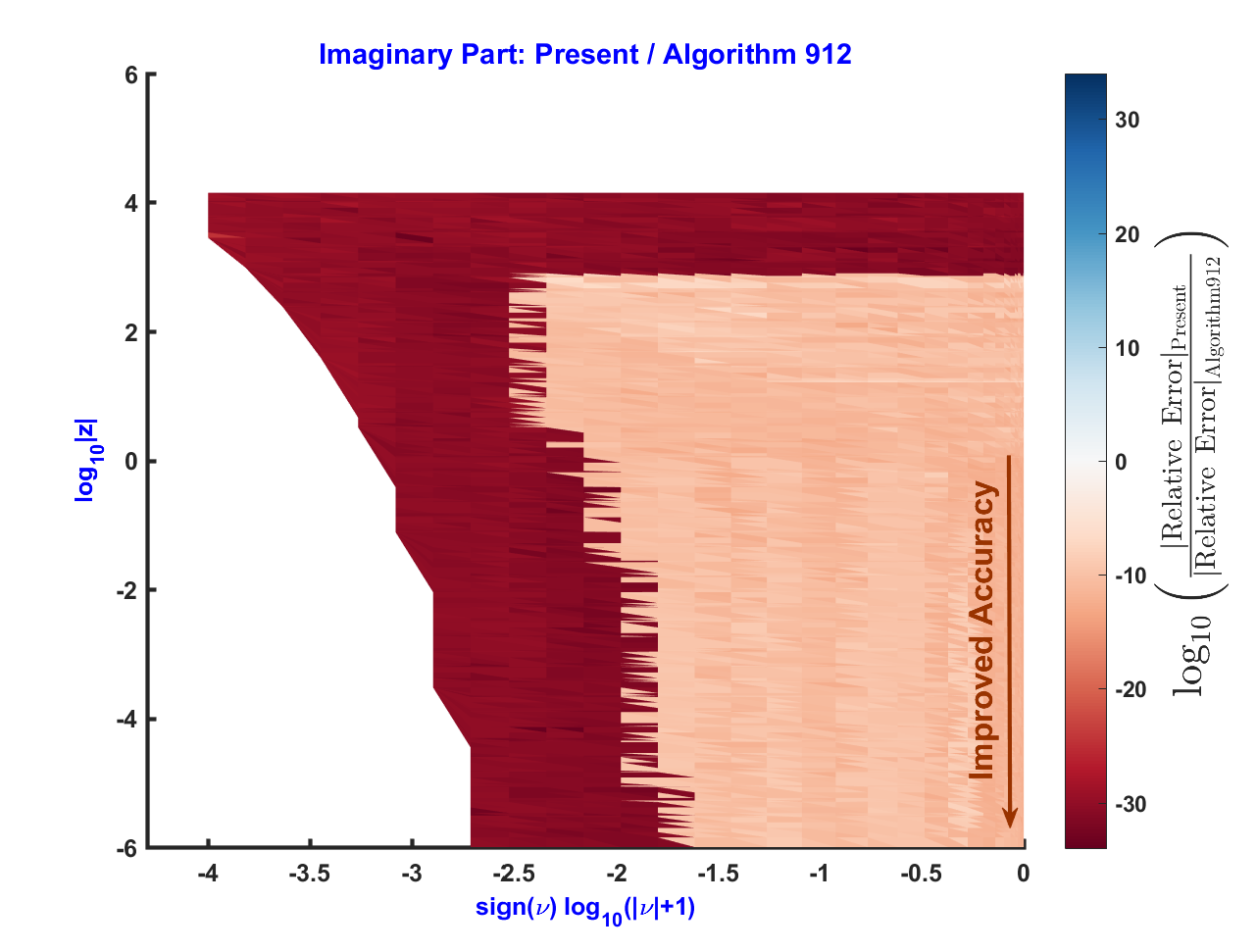}}
   \caption{Colormap plots of the base-10 logarithm of the \textit{ratio} of componentwise
relative errors (present/Algorithm~912) in computing the real part~(a) and
imaginary part~(b) of \(J_{-|\nu|}(z)\) in quadruple precision, using
\texttt{Maple} values as the reference.}
    \label{fig:qp_j_nu_im_912}
\end{figure}

\noindent
The comparison with Algorithm~912 in quadruple precision reveals a more pronounced distinction than in double precision. As illustrated in Figures~\ref{fig:qp_j_nu_re_912} and~\ref{fig:qp_j_nu_im_912}, the present implementation consistently achieves significantly smaller relative errors than Algorithm~912, even within its nominal range of applicability. This difference reflects the enhanced numerical stability and adaptive selection of representations in the present method. Moreover, while Algorithm~912 is confined to a restricted region of the \((\nu,z)\)-plane, the present algorithm continues to provide stable and accurate evaluations across the substantially enlarged quadruple-precision test domain.

\FloatBarrier
\subsubsection{Bessel Functions of the Second Kind \texorpdfstring{$\Ynu(z)$}{Ynu(z)}} \label{sec:acc_qp_Y}

\begin{figure}[htbp]
    \centering
    \subfloat[]{\includegraphics[width=0.48\linewidth]{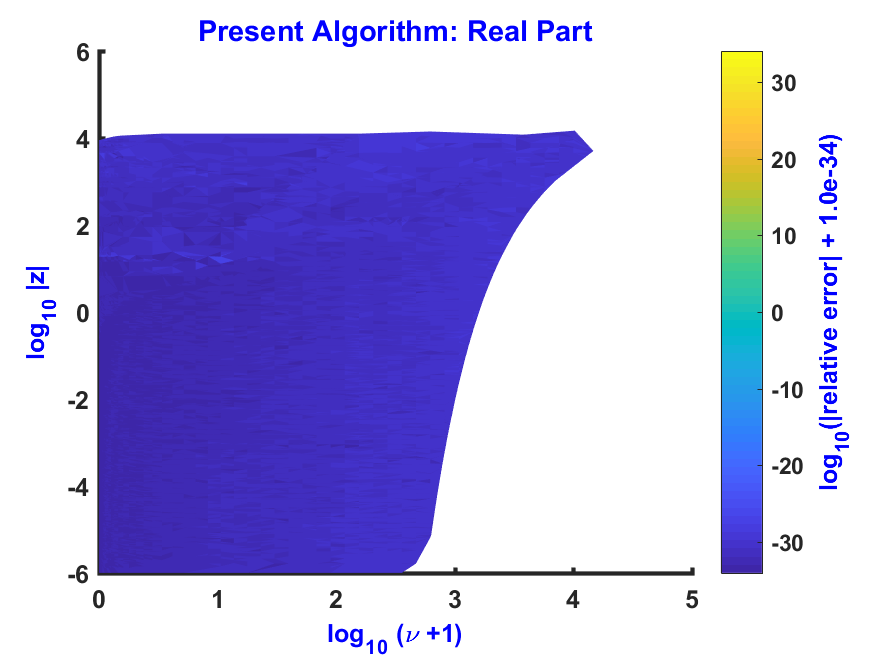}}
    \subfloat[]{\includegraphics[width=0.48\linewidth]{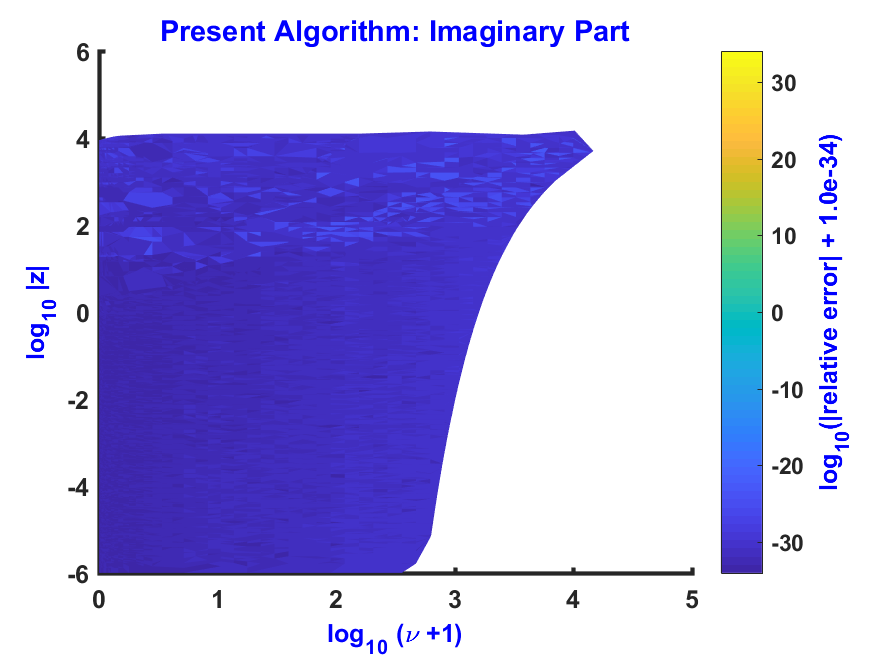}}
    \caption{Colormap plots of the base-10 logarithm of the relative error in computing the real
    part~(a) and imaginary part~(b) of \(\Ynu(z)\) in quadruple precision,
    using \texttt{Maple} values as the reference.}
    \label{fig:Acc_rel_err_Y}
\end{figure}

The quadruple-precision implementation of \(\Ynu(z)\) also benefits from the extended dynamic range and enhanced precision. Over a dense grid of 278{,}460 test points in the enlarged \(\nu\)–\(|z|\) domain, the method remains stable and attains very small relative errors for positive orders.
Figure~\ref{fig:Acc_rel_err_Y} shows the relative error in computing the real part~(a) and imaginary part~(b) of \(\Ynu(z)\) for positive orders, with \texttt{Maple} values used as the reference.

\begin{figure}[htbp]
    \centering
    \subfloat[]{\includegraphics[width=0.48\linewidth]{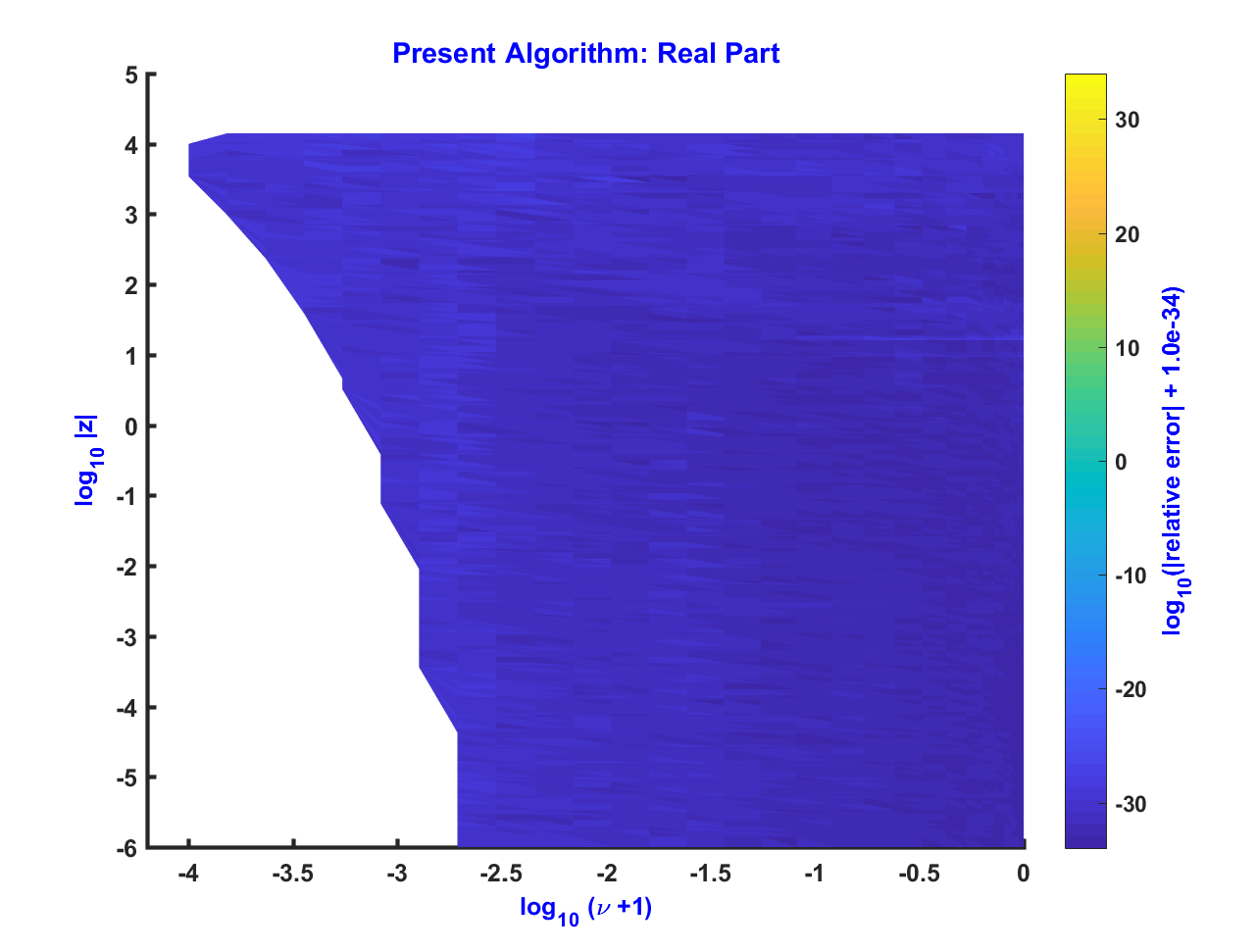}}
    \subfloat[]{\includegraphics[width=0.48\linewidth]{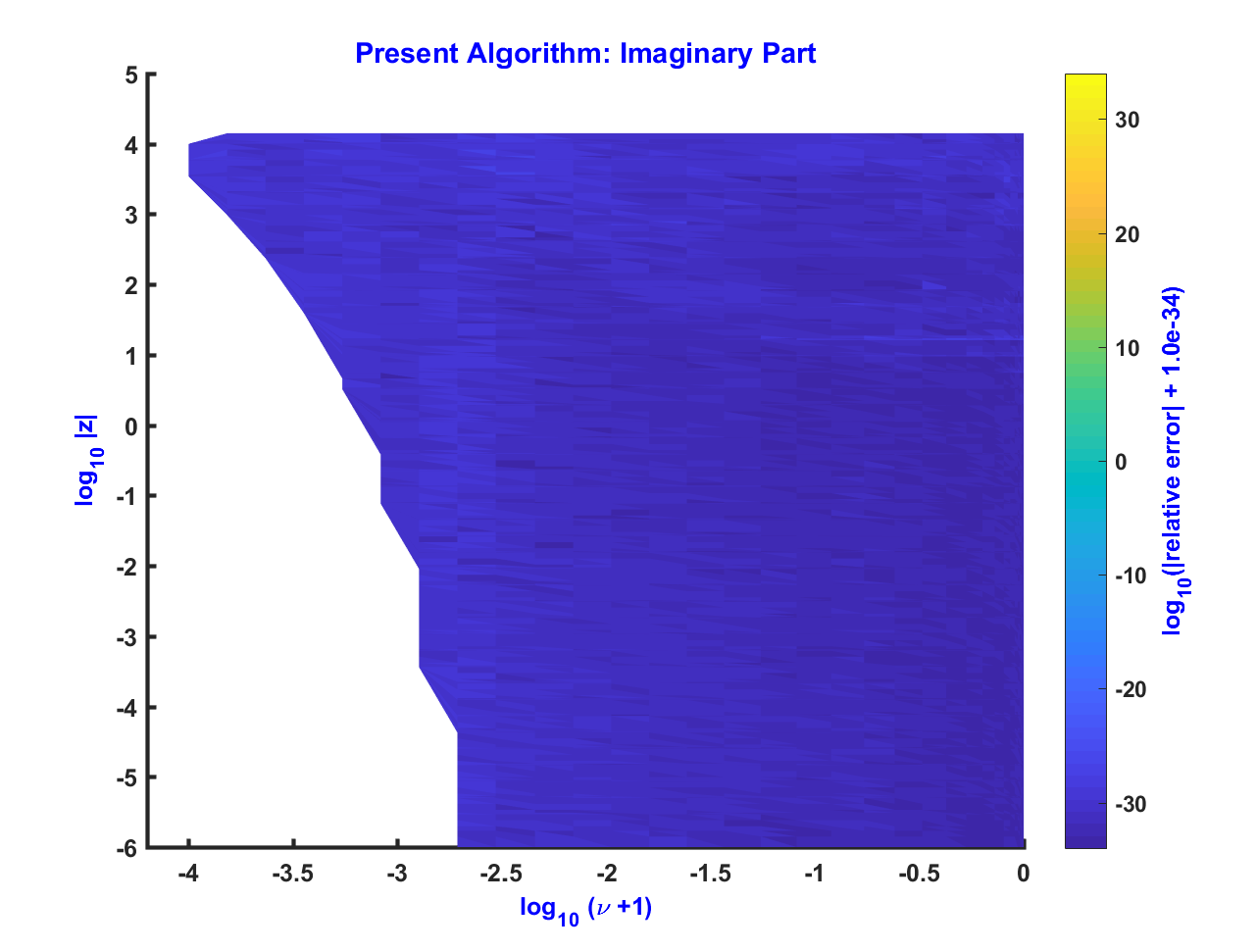}}
    \caption{Colormap plots of the base-10 logarithm of the relative error in computing the real
    part~(a) and imaginary part~(b) of \(Y_{-|\nu|}(z)\) in quadruple precision,
    using \texttt{Maple} values as the reference.}
    \label{fig:Acc_rel_err_Y_neg_nu}
\end{figure}

For negative orders, a test grid of 277{,}975 points, primarily distributed on a logarithmic scale, is used to assess the accuracy of computing \(Y_{-|\nu|}(z)\) in quadruple precision with the present implementation. Figure~\ref{fig:Acc_rel_err_Y_neg_nu} shows that the algorithm maintains high accuracy throughout the enlarged domain, using high-precision \texttt{Maple} values as the reference, and exhibits no loss of numerical stability in the negative-order region.

\begin{figure}[htbp]
    \centering
    \subfloat[]{\includegraphics[width=0.48\linewidth]{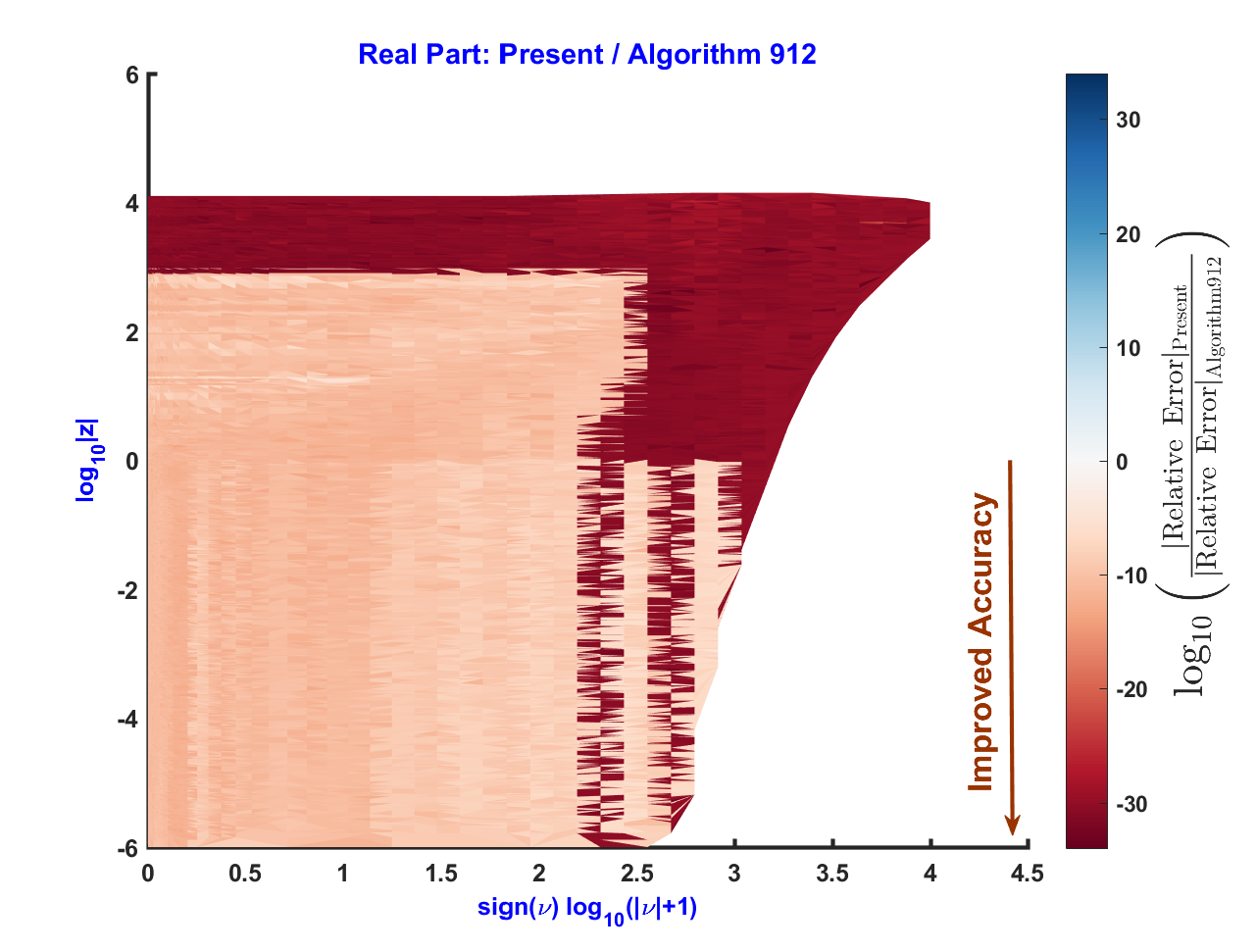}}
    \subfloat[]{\includegraphics[width=0.48\linewidth]{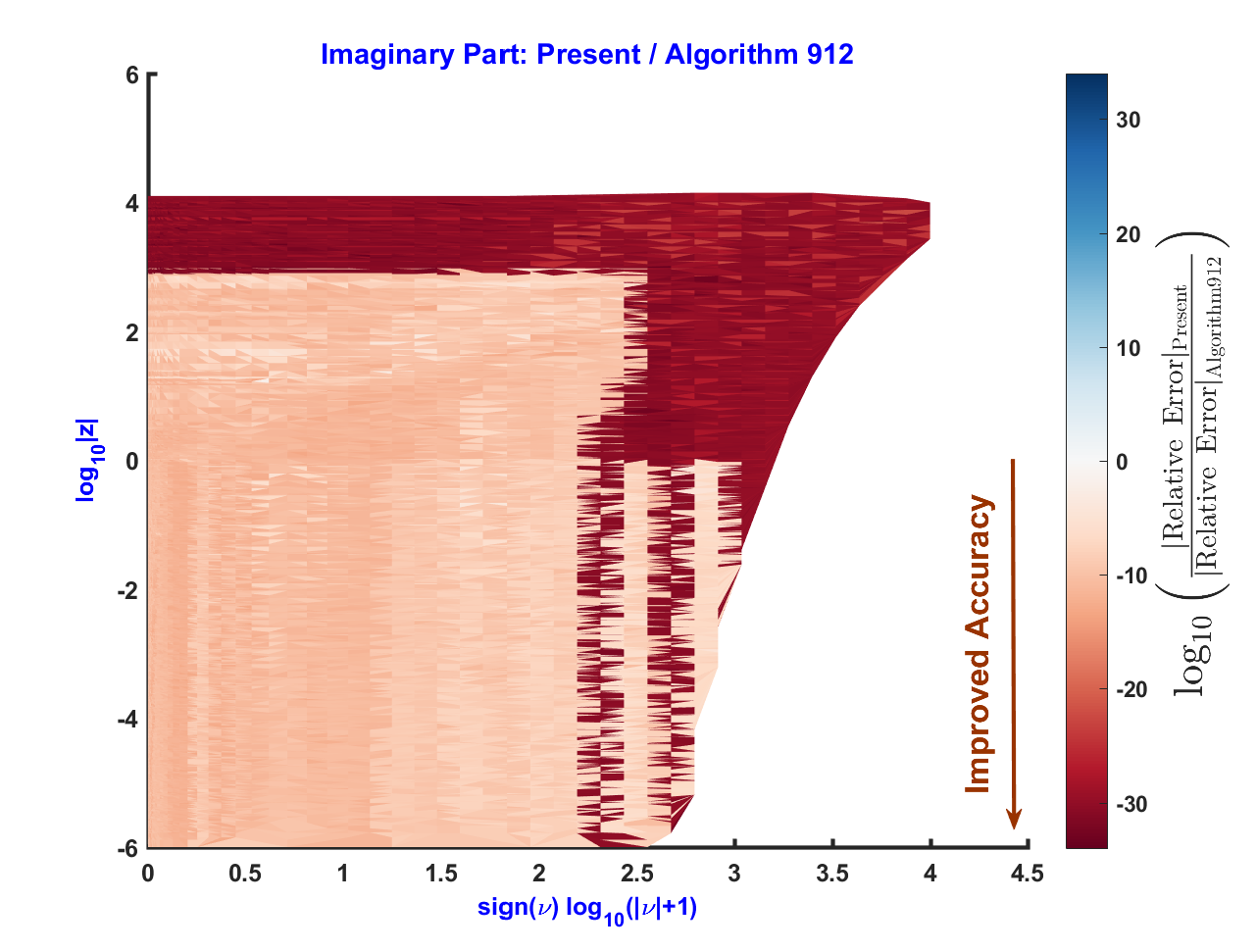}}
   \caption{Colormap plots of the base-10 logarithm of the \textit{ratio} of componentwise
relative errors (present/Algorithm~912) in computing the real part~(a) and
imaginary part~(b) of \(Y_{\nu}(z)\) in quadruple precision, using
\texttt{Maple} values as the reference.}
    \label{fig:Acc_Y_nu_912_qp}
\end{figure}

\begin{figure}[htbp]
    \centering
    \subfloat[]{\includegraphics[width=0.48\linewidth]{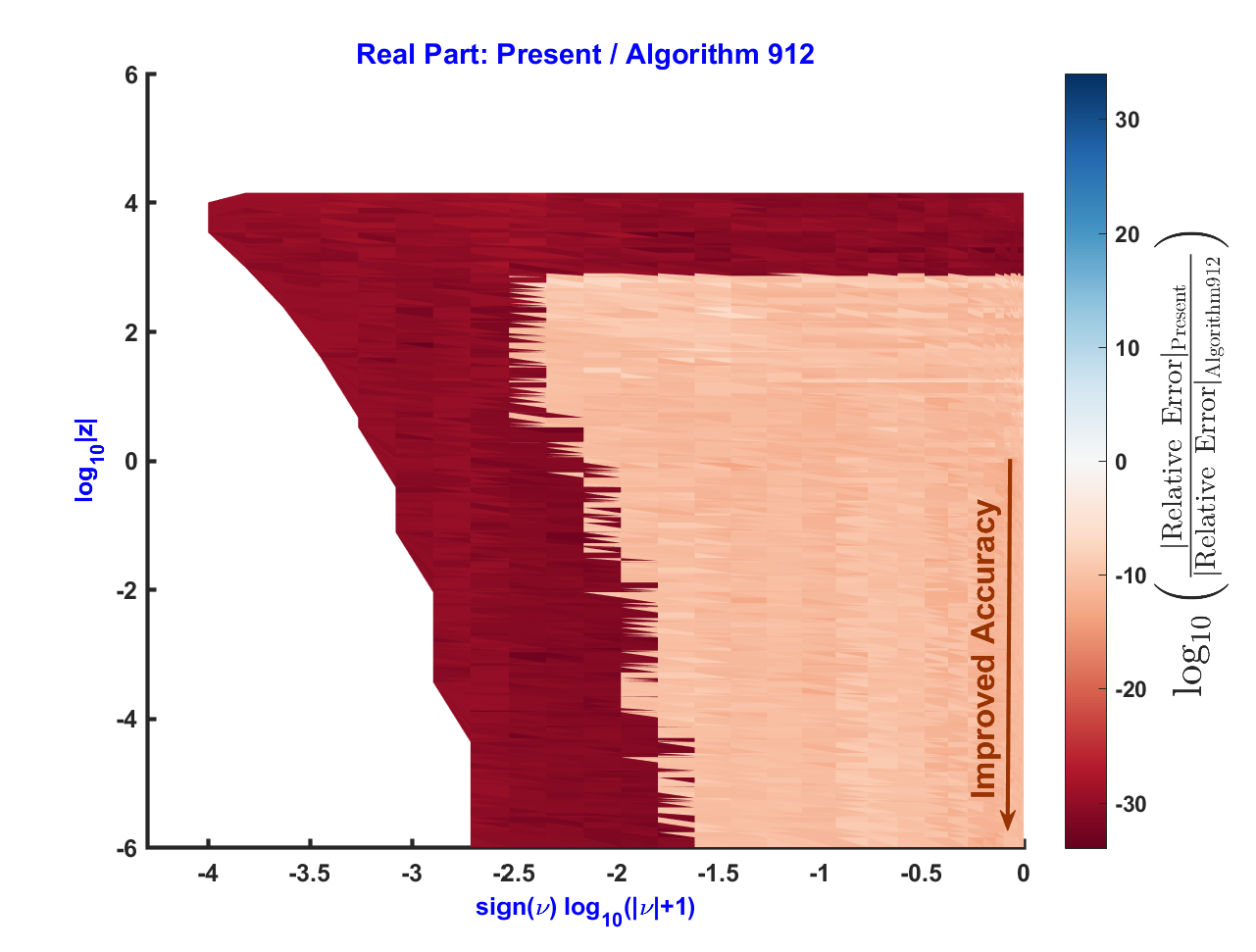}}
    \subfloat[]{\includegraphics[width=0.48\linewidth]{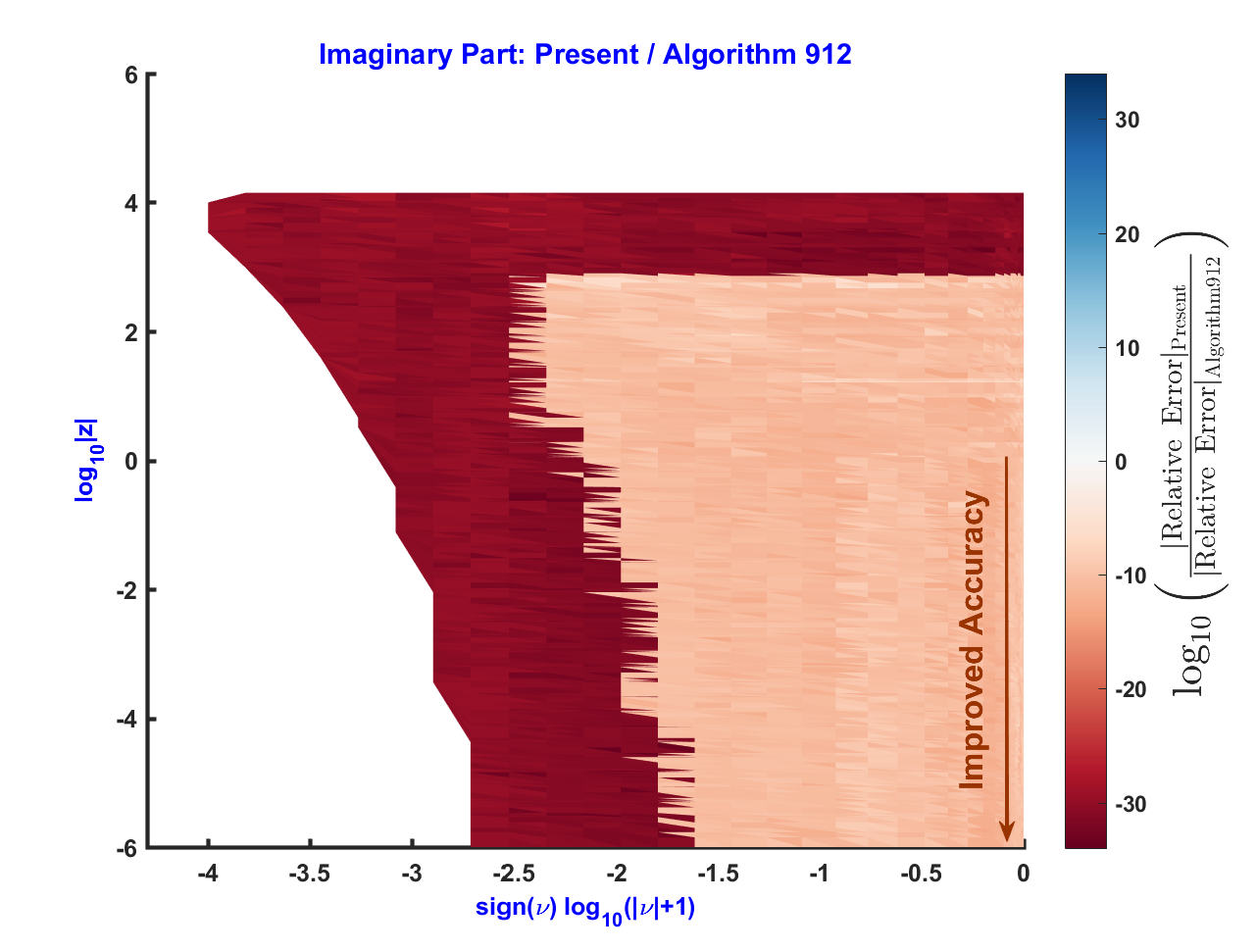}}
   \caption{Colormap plots of the base-10 logarithm of the \textit{ratio} of componentwise
relative errors (present/Algorithm~912) in computing the real part~(a) and
imaginary part~(b) of \(Y_{-|\nu|}(z)\) in quadruple precision, using
\texttt{Maple} values as the reference.}
    \label{fig:Acc_Y_neg_nu_912_qp}
\end{figure}

Consistent with the behavior observed for \(\Jnu(z)\), comparisons with Algorithm~912, as shown in Figures~\ref{fig:Acc_Y_nu_912_qp} and~\ref{fig:Acc_Y_neg_nu_912_qp}, reveal the same qualitative trends: the present implementation achieves significantly higher accuracy, while Algorithm~912 exhibits both reduced accuracy and a more limited domain of applicability.

\FloatBarrier
\section{Efficiency Benchmarking}

\subsection{Comparison with Algorithm~644}
To assess the computational efficiency of the present algorithms for \(\Jnu(z)\) and \(\Ynu(z)\), we employ the same benchmarking methodology 
used previously for the modified Bessel functions \(I_{\nu}(z)\) and \(K_{\nu}(z)\) in Parts~I and~II of this series. In particular, a systematic 
performance comparison is carried out against Algorithm~644 using the extensive datasets noted in Section~\ref{sec:acc_dp}, ensuring a thorough 
and consistent evaluation across all members of the Bessel family considered in this work.

Execution times for both algorithms are measured under identical computational conditions, enabling a direct and fair comparison of performance. Again, since 
Algorithm~644 is restricted to double-precision arithmetic and offers no support for quadruple precision, all efficiency measurements are carried out in double precision. The resulting benchmarks emphasize the computational advantages of the present algorithms, including their reduced execution time and their ability to operate reliably in parameter regimes where Algorithm~644 either fails or exhibits degraded numerical performance.

As in the accuracy verification studies, efficiency tests are conducted using the same set of test points employed for double-precision. Each benchmark 
consists of evaluating the complete dataset 50 times, with the entire process repeated 21 times. To mitigate the effects of system-level variability and timing noise, the minimum execution time among the 21 repetitions is reported, following standard best practices for performance benchmarking \cite{BenchmarkTools}.  Timing measurements are obtained using the \texttt{SYSTEM\_CLOCK} intrinsic, which maps to the high-resolution \texttt{QueryPerformanceCounter} system call in GNU~Fortran, thereby providing sufficient resolution for accurate measurement of the short runtimes 
associated with pointwise function evaluations.

All benchmarking experiments are performed under controlled conditions on a fixed hardware platform, with background system activity minimized to reduce 
external sources of timing variability. Performance measurements are obtained using the following \texttt{Fortran} compilers:

\begin{itemize}
    \item GNU Fortran (i686-posix-dwarf-rev0, Built by MinGW-W64 project) 8.1.0,
    \item GNU Fortran (Rev3, Built by MSYS2 Project) 12.1.0,
    \item NAG Fortran Compiler Release 7.1 (Hanzomon) Build 7110,
    \item Intel(R) Fortran Intel(R) 64 Version 2021.9.0 Build 20230302\_000000,
    \item \texttt{IFX} (LLVM-based), from Intel(R) Fortran 64 Version 2021.9.0.
\end{itemize}

All test cases are compiled using a range of standard compiler optimization levels (\texttt{-O0}, \texttt{-O1}, \texttt{-O2}, and \texttt{-O3}) in order to assess the impact of compiler optimizations on performance. Benchmark results are reported with an emphasis on reproducibility, and all relevant system and compiler specifications are documented accordingly.

Since \(\Jnu(z)\) and \(\Ynu(z)\) are computed primarily through the evaluation of \(I_{\nu}(z)\) and \(K_{\nu}(z)\), execution-time measurements are  reported only for the entire datasets. The summarized results, presented in Table~\ref{tab:execution_time_comparison_transposed}, show that the present algorithm consistently outperforms Algorithm~644 over the full tested domain, independent of 
the compiler used and the selected optimization level. When expressed as a percentage of the computational time required by Algorithm~644, the present implementation achieves execution times ranging from approximately 35\% to 67\% for \(\Jnu(z)\) and from approximately 44\% to 72\% for \(\Ynu(z)\), depending on the compiler and optimization settings.

\begin{figure}[htbp]
  \centering
  \subfloat[]{%
    \includegraphics[width=0.47\linewidth]{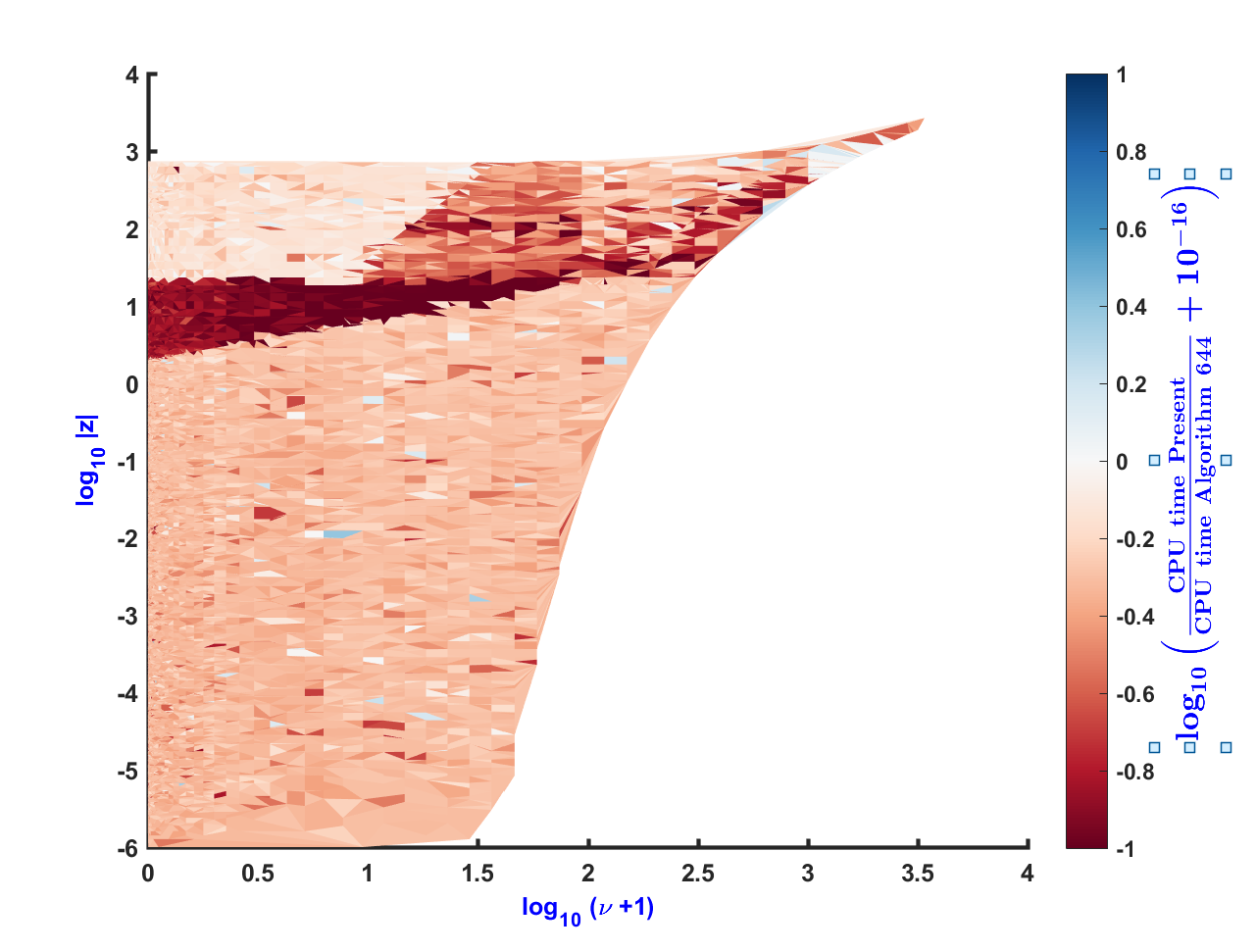}%
    \label{fig:pt_by_pt_timing_J}%
  }\hfill
  \subfloat[]{%
    \includegraphics[width=0.47\linewidth]{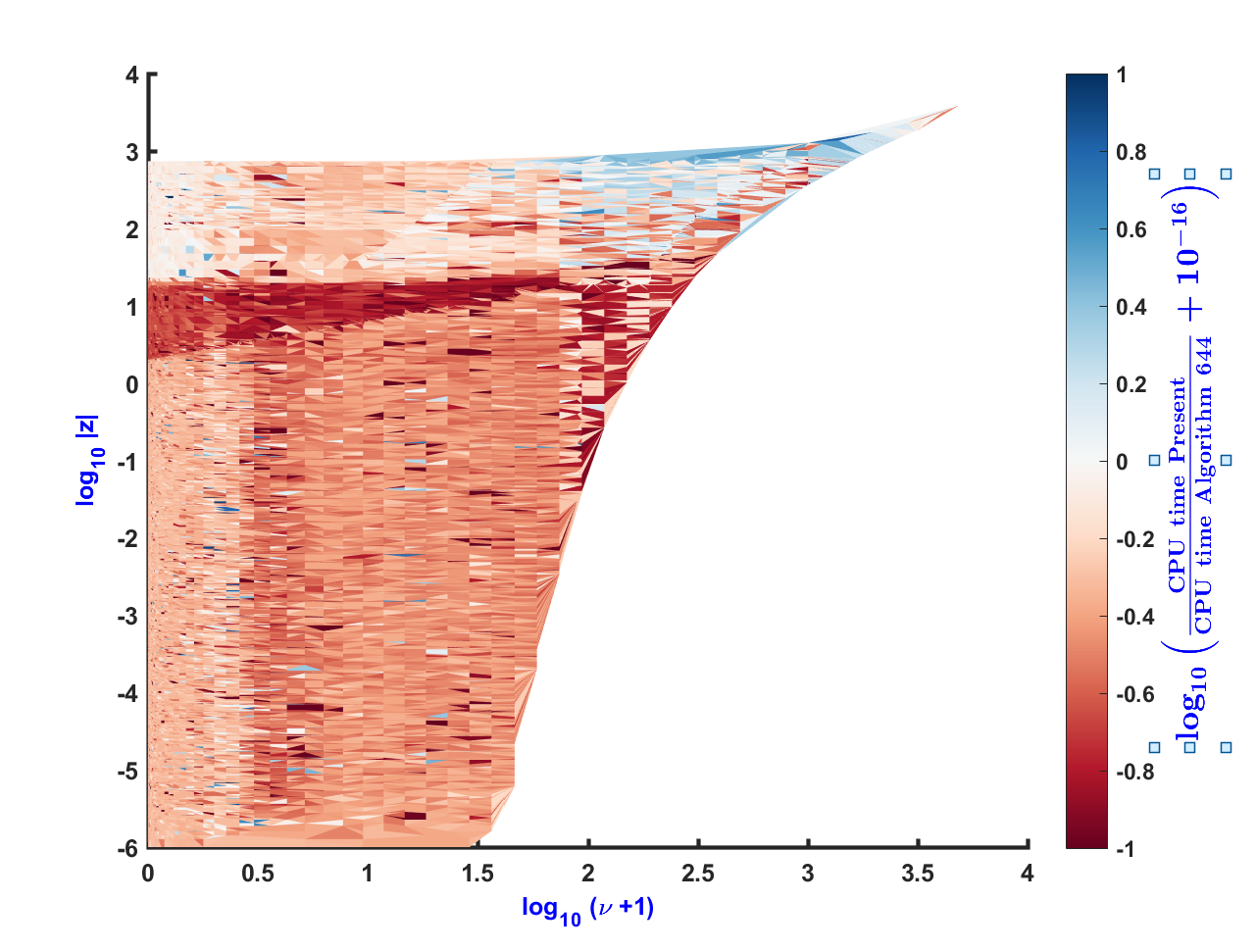}%
    \label{fig:pt_by_pt_timing_Y}%
  }
 \caption{Two-dimensional colormap plots of the base-10 logarithm of the ratio of
the execution time of the present algorithm to that of Algorithm~644 across the
entire computational domain for: \textbf{(a)} \(\Jnu(z)\) and \textbf{(b)}
\(\Ynu(z)\). Computations were performed on an
Intel\textsuperscript{\textregistered} Core\textsuperscript{TM} i7-6600U CPU
(2.60--2.81\,GHz) using GNU Fortran (\textit{gfortran},
i686-posix-dwarf-rev0, built by the MinGW-W64 Project), version~8.1.0.}

\end{figure}

As in the previous two parts of this series, additional insight into the performance characteristics of the algorithms is obtained by examining the
computational cost on a pointwise basis, rather than relying solely on averages over the full domain. Figure~\ref{fig:pt_by_pt_timing_J}
presents a colormap of the base-10 logarithm of the \emph{ratio} of execution times for a single evaluation of \(\Jnu(z)\) using the present algorithm
relative to Algorithm~644. For each point, the reported value corresponds to the minimum execution time observed across multiple independent runs for each
algorithm, thereby reducing the influence of transient system-level fluctuations. Figure~\ref{fig:pt_by_pt_timing_Y} shows the corresponding
execution-time ratio for \(\Ynu(z)\).

Although the present algorithm demonstrates superior performance over the majority of the domain, there exist small, isolated regions where Algorithm~644 attains shorter execution times. Some of these regions are located near the underflow/overflow boundary, where Algorithm~644 effectively avoids computation by returning immediate zero values. Other isolated regions are associated primarily with extremely large orders and arguments, suggesting opportunities for further performance optimization in these extreme parameter regimes.

As illustrated in Figure~\ref{fig:pt_by_pt_timing_Y}, the performance gains of the present algorithm are most pronounced in the intermediate-\(z\) and intermediate-\(\nu\) regimes. In these regions, the algorithm for \(\Jnu(z)\) and consequently \(\Ynu(z)\) exploits  either directly or indirectly the efficient small-\(z\) series expansions and direct absolute backward recurrence relations, respectively. 

The average execution time per evaluation in quadruple precision is roughly an order of magnitude slower than in double precision.

\begin{table}[ht]
\centering
\caption{Execution times of the present algorithms relative to Algorithm~644 (in percent) for \(\Jnu(z)\) and \(\Ynu(z)\) using double-precision
arithmetic. Results are shown for different \texttt{Fortran} compilers and optimization levels.}
\label{tab:execution_time_comparison_transposed}

\renewcommand{\arraystretch}{1.2}
\begin{tabular}{lccccc} 
\hline
\textbf{Compiler } &
\textbf{Gfortran 8.1} &
\textbf{Gfortran 12.1} &
\textbf{NAG 7.1} &
\textbf{IFORT 2021} &
\textbf{IFX 2021} \\
\hline

\multicolumn{6}{l}{\textbf{$\Jnu(z)$}} \\
\hline
\texttt{-O3} & 35.1\% & 60.9\% & 57.7\% & 58.9\% & 55.7\% \\
\texttt{-O2} & 35.2\% & 60.9\% & 58.7\% & 58.4\% & 55.7\% \\
\texttt{-O1} & 34.6\% & 60.7\% & 56.7\% & 66.2\% & 55.0\% \\
\texttt{-O0} & 35.9\% & 53.9\% & 66.5\% & 49.8\% & 55.4\% \\

\hline
\multicolumn{6}{l}{\textbf{$\Ynu(z)$}} \\
\hline
\texttt{-O3} & 45.3\% & 73.2\% & 65.8\% & 65.9\% & 66.4\%\\
\texttt{-O2} & 45.0\% & 71.5\% & 63.2\% & 68.1\% & 66.6\%\\
\texttt{-O1} & 44.1\% & 70.3\% & 64.4\% & 65.1\% & 65.6\%\\
\texttt{-O0} & 46.9\% & 63.5\% & 55.1\% & 71.8\% & 66.1\% \\
\hline
\end{tabular}
\end{table}

\subsection{Comparison with Algorithm~912}

Table~\ref{tab:timing_ratio_912} summarizes the total execution times for the
present implementation and Algorithm~912. The ratios show that the present
algorithm requires only a small fraction of the time required by Algorithm~912,
with speedups ranging from approximately \(22\times\) to \(256\times\) over the
tested cases.

\begin{table}[ht]
\centering
\caption{Execution times, ratios, and speedups relative to Algorithm~912. The
ratio is \(T_{\rm present}/T_{912}\), and the speedup is \(T_{912}/T_{\rm present}\).}
\label{tab:timing_ratio_912}

\renewcommand{\arraystretch}{1.2}
\begin{tabular}{lcccc}
\hline
\textbf{Function \& precision} &
\textbf{Present (s)} &
\textbf{Algorithm~912 (s)} &
\textbf{\(T_{\rm present}/T_{912}\)} &
\textbf{Speedup} \\
\hline
\(J\), double precision (\(+\nu\))       & 1.683  & 430.84 & \(3.91\times10^{-3}\) & \(256.0\times\) \\
\(J\), double precision (\(-\nu\))       & 2.48   & 565    & \(4.39\times10^{-3}\) & \(227.8\times\) \\
\(Y\), double precision (\(+\nu\))       & 3.07   & 455    & \(6.75\times10^{-3}\) & \(148.2\times\) \\
\(Y\), double precision (\(-\nu\))       & 3.38   & 354    & \(9.55\times10^{-3}\) & \(104.7\times\) \\
\hline
\(J\), quadruple precision (\(+\nu\))    & 18.55  & 1399.7 & \(1.33\times10^{-2}\) & \(75.5\times\) \\
\(J\), quadruple precision (\(-\nu\))    & 74.7   & 2500   & \(2.99\times10^{-2}\) & \(33.5\times\) \\
\(Y\), quadruple precision (\(+\nu\))    & 97.35  & 2658.4 & \(3.66\times10^{-2}\) & \(27.3\times\) \\
\(Y\), quadruple precision (\(-\nu\))    & 115.8  & 2574.8 & \(4.50\times10^{-2}\) & \(22.2\times\) \\
\hline
\end{tabular}
\end{table}

\begin{figure}[htbp]
  \centering
  \subfloat[]{%
    \includegraphics[width=0.47\linewidth]{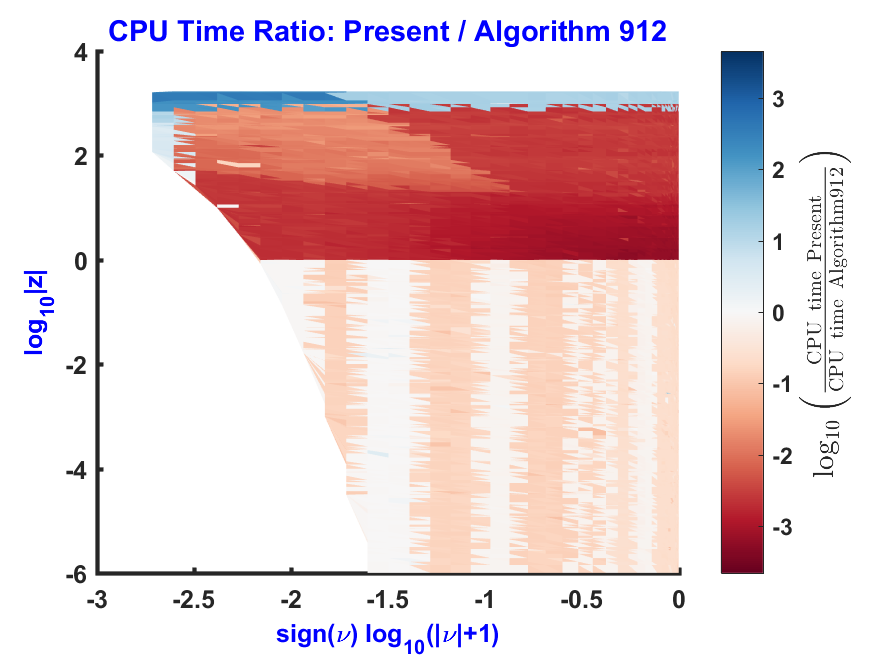}%
   \label{fig:pt_by_pt_timing_J_neg_nu_912}%
  }\hfill
  \subfloat[]{%
    \includegraphics[width=0.47\linewidth]{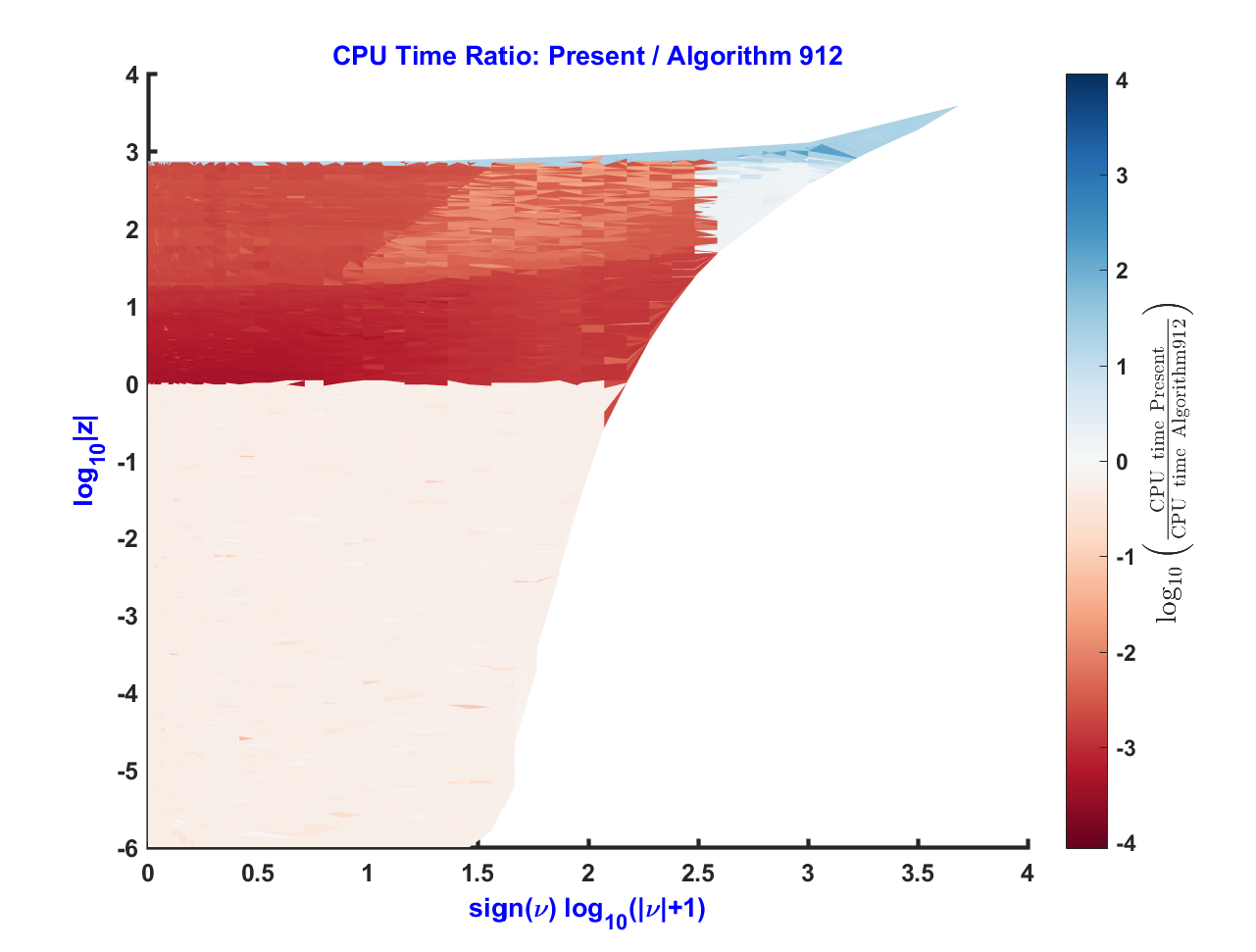}%
    \label{fig:pt_by_pt_timing_J_912}%
  }
\caption{Two-dimensional colormap plots of the base-10 logarithm of the
\textit{ratio} of the execution time of the present algorithm to that of
Algorithm~912 across the computational domain for: \textbf{(a)}
\(J_{-|\nu|}(z)\) and \textbf{(b)} \(J_{\nu}(z)\). Computations were performed
on an Intel\textsuperscript{\textregistered} Core\textsuperscript{TM} i7-6600U
CPU (2.60--2.81\,GHz) using GNU Fortran (\textit{gfortran},
i686-posix-dwarf-rev0, built by the MinGW-W64 Project), version~8.1.0.}
\end{figure}

Figures~\ref{fig:pt_by_pt_timing_J_neg_nu_912} and
\ref{fig:pt_by_pt_timing_J_912} show representative double-precision
pointwise timing comparisons for \(J_{-|\nu|}(z)\) and \(J_{\nu}(z)\),
respectively. The pointwise timing results show that Algorithm~912 becomes
substantially more expensive in regions of moderate and large \(|z|\). This
behavior is consistent with the structure of Algorithm~912: outside the
small-argument region, the method relies on numerical integration, including
the construction of steepest-descent contours, root-finding iterations, and
adaptive quadrature. These steps introduce substantial overhead for each
function evaluation.

The isolated regions where Algorithm~912 appears comparatively faster should
not be interpreted as genuine performance advantages. In many such cases, the
routine terminates early because its internal range or reliability checks are
triggered. Thus, the shorter execution time corresponds to an incomplete or
rejected evaluation rather than to a successful computation of the requested
Bessel function.

Consequently, the timing results must be interpreted together with the accuracy and applicability results. In regions where Algorithm~912 returns valid outputs,
it is consistently more expensive than the present implementation; outside those regions, the present algorithm continues to provide accurate and stable values
over the tested IEEE double- and quadruple-precision domains. Generating analogous full-domain pointwise timing maps for all \(Y_\nu\) cases and for
quadruple precision was not practical because of the prohibitive computational cost of Algorithm~912, especially in regions where its numerical-integration branch is activated.

\section{Conclusions}
\label{sec:conclusions}

This work completes a unified effort aimed at the efficient multiprecision computation of Bessel functions. In this final part, we have presented robust,
highly accurate, and computationally efficient algorithms, together with \texttt{Fortran} implementations, for evaluating the regular Bessel functions
\(\Jnu(z)\) and \(\Ynu(z)\) for complex arguments and for both positive and \textit{negative} real orders. 

The proposed methods overcome several limitations of existing approaches, most notably Algorithm~644, by substantially extending the computational domain into parameter regimes where that method fails or exhibits degraded numerical performance. Extensive verification against high-precision \texttt{Maple}
reference values demonstrates high accuracy across the tested domain for all points with representable reference values, while comprehensive benchmarking confirms consistent performance advantages over Algorithm~644 across all tested compilers and optimization levels.

\medskip
\noindent
Comparisons with Algorithm~912 further highlight the advantages of the present approach. Although Algorithm~912 provides extended functionality, including support for complex orders and quadruple-precision, its applicability is restricted to limited regions of the \((\Re \nu,z)\)-plane and its computational cost increases substantially in challenging parameter regimes. In particular, for large values of \(|z|\), the reliance of Algorithm~912 on numerical integration and auxiliary root-finding procedures leads to significantly increased computational effort.

The availability of native quadruple-precision support addresses the growing need for enhanced numerical stability and extended dynamic range in demanding
scientific applications. Moreover, the modular structure of the algorithms facilitates straightforward adaptation to other numerical formats beyond double and quadruple
precision with minimal implementation effort.

By providing accurate, robust, and efficient algorithms for \(\Jnu(z)\) and \(\Ynu(z)\), this work completes a unified, multiprecision computational
framework for Bessel and modified Bessel functions of real order and complex argument. The resulting software package will be released under an open-source
license to support reproducible research and to enable further developments in high-accuracy special-function computation.

\section*{Software and reproducibility}

The complete \texttt{Fortran} source code, test data, plotting scripts, compiler
invocation lines, and a \texttt{Makefile} will be made available at
\texttt{https://github.com/mofrehzaghloul/}. The repository will also include
instructions for reproducing the accuracy tests, timing experiments, and figures
reported in this paper.

\section*{Acknowledgments}
This work was supported by the UAE University UPAR research grant \textbf{G00004995 (2024)}.

The authors gratefully acknowledge valuable comments and suggestions from Miguel Hermanns (Departamento de Mecánica de Fluidos y Propulsión Aeroespacial, Madrid, Spain).

The first author also gratefully acknowledges the hospitality and support of S.~G.~Johnson and the Department of Mathematics at the Massachusetts Institute of Technology during his visits as a Visiting Professor in Spring 2025 and Summer 2026.



\bibliographystyle{siamplain}
\bibliography{references}

\end{document}